\documentclass[12pt,a4paper]{article}

\usepackage{pdflscape}
\usepackage{amsmath}
\usepackage{amssymb}
\usepackage{latexsym}
\usepackage{srcltx}
\usepackage{graphics}
\usepackage{color}
\usepackage{epsfig}
\usepackage{color}
\usepackage[ruled,vlined]{algorithm2e}

\newcommand{\re}{{\mathbb R}}
\newcommand{\n}{{\mathbb N}}

\newcommand{\z}{{\mathbb Z}}

\newcommand{\bx}{{\boldsymbol x}}

\newcommand{\ba}{{\boldsymbol a}}

\newcommand{\bz}{{\boldsymbol z}}
\newcommand{\be}{{\boldsymbol e}}

\newcommand{\bv}{{\boldsymbol v}}
\newcommand{\bu}{{\boldsymbol u}}

\newcommand{\bell}{{\boldsymbol \ell}}
\newcommand{\bzero}{{\bf 0}}

\newcommand{\wA}{\widetilde A}
\newcommand{\wX}{\widetilde X}

\newtheorem{theorem}{Theorem}
\newtheorem{prop}{Proposition}
\newtheorem{lemma}{Lemma}
\newtheorem{cor}{Corollary}
\newtheorem{remark}{Remark}
\newtheorem{ex}{Example}
\newtheorem{defi}{Definition}

\date{}

\author{Nicola Guglielmi
\thanks{Dipartimento di Ingegneria Scienze Informatiche e Matematica,
University of L'Aquila, Italy; Gran Sasso Science Institute, via Crispi 7,
L'Aquila, Italy, {e-mail: \tt\small guglielm@univaq.it }}
 and
Vladimir Yu.~Protasov
\thanks{Dipartimento di Ingegneria Scienze Informatiche e Matematica,
University of L'Aquila, Italy; Department
of Computer Science of Higher School of Economics, Moscow, Russia,  {e-mail: \tt\small
v-protassov@yandex.ru}}}

\title{On the closest stable/unstable nonnegative matrix \\
and related stability radii
\thanks{
 The second author is supported by RSF Grant 17-11-01927 }}

\begin{document}
\maketitle

\begin{abstract}

We consider the problem of computing the closest 
stable/unstable non-negative matrix to a given real matrix. This problem is important  in the study of linear dynamical systems, numerical methods, etc.  The distance between matrices is measured in the Frobenius norm. The problem is addressed for two types of stability: the Schur stability (the matrix is stable if its spectral radius is smaller than one) and the Hurwitz stability (the matrix is stable if its spectral abscissa is negative). We show that the closest unstable matrix can always be explicitly found.  For the closest stable matrix, we present an iterative algorithm
which  converges to a local minimum with a linear rate. It is shown  that 
the total number of local  minima can be exponential in the dimension. Numerical results and the complexity estimates are presented.

\smallskip

\smallskip

\begin{flushright}
\noindent  \textbf{AMS 2010} {\em subject
classification: 15B48, 34K20, 90C26, 65F15} 
\end{flushright}

\end{abstract}
\bigskip

\begin{center}
\textbf{1. Introduction}
\end{center}
\bigskip

 The problems of finding the closest
stable matrix ({\em stabilizing problem}) or the closest unstable matrix 
({\em destabilizing problem}) are very important in many applications such as 
the
analysis of differential equations, linear dynamical
systems, electrodynamics, etc. In this paper we focus on those problems in the set of non-negative matrices and call them {\em positive stabilizing/destabilizing problems}.  They are needed in the study of positive linear systems which are  widely applied in the multiagent problems, population dynamics (matrix population models), mathematical economics (Leontief model), etc. 
 Non-negativity ensures certain advantages for this problem such as the special spectral properties
of matrices guaranteed by the Perron-Frobenius theory. On the other hand,  
it brings extra constraints (the non-negativity of $d^2$ entries, where $d$ is the dimension) which complicates the problem a lot. 

For the Schur stability, the problem of non-negative stabilizing of a given matrix $A$ 
consists in finding an entriwise non-negative matrix $X$ such such that 
$\rho(X) \le 1$ and the distance to the matrix $\|X-A\|$ is minimal. 
As usual, $\rho(X)$ denotes the spectral radius of the matrix, which is 
the maximum modulus of its eigenvalues.  As a rule, such stabilizing problems are  notoriously
hard due to properties of the spectral radius~$\rho(X)$ as a
function of the matrix~$X$. The function~$\rho(X)$ is neither convex nor concave, it
may be non-differentiable and even non-Lipschitz  at some points. 
This makes all methods of convex or smooth optimization hardly applicable. Basically, even 
finding a locally closest stable matrix is hard. 
In many situations we cannot hope for the global optimality due to 
a large number of local minima. Methods of matrix stabilization
(without non-negativity assumption) were presented in~\cite{A, By, GM, GS,  MMS, NP, ONVD}.

In the stabilizing problem a lot depends on the norm we measure the distance 
$\|X-A\|$. Usually it is either Euclidean or Frobenius norm; the problem is hard in both 
those norms~\cite{By, GS, ONVD}.  
In the recent paper~\cite{NP} it has been shown that in the $L_{\infty}$ matrix norm 
(equal to the largest $L_1$ norm of  rows of the matrix)
the problems of positive stabilizing/distabilizing  both have  surprisingly simple  
solutions and there exist  efficient algorithms that find their global minimuma. 
In applications, however, the $L_{\infty}$ matrix norm has some disadvantages:
it is non-smooth,  badly correlated with the Euclidean norm, etc.  
That is why, many researches prefer the Frobenius norm, which is the sum of squares 
of the matrix components. Actually, the Frobelius norm is merely the vector Euclidean norm 
in the $d^2$ dimensional space of matrices.  In this paper we deal with the positive 
stabilizing and destabilizing problems in the Frobenius norm. 
We show how to find explicitly the closest non-negative unstable matrix. 
In fact, it was observed in the literature that the destabilizing problem is usually 
simpler than the stabilizing one. In some favorable cases our method can be extended 
for finding the closest stable matrix as well, but in general we may hope on finding 
local minima only. To this end,  we develop an iterative relaxation scheme that converges to a 
local minimum. It is computationally simple and uses only a standard quadratic programming routine. This makes the method applicable even in high dimensions. Another advantage is the stability of the method with respect to matrices: the 
algorithm works 
equally well even if the matrix in a current iteration has a multiple leading eigenvalue (which happens often) and its spectral radius is non-Lipscitz. 
In  practice, the new method  converges extremely fast, which is demonstrated in  numerical examples.  We prove that the rate of convergence is always linear and, moreover, 
if the limit matrix is strictly positive then it gives a global minimum. In general the limit matrix may have zero components, in which case it gives only a local minimum in general. In this case, a 
question arises about the possible number of local minima. 
We construct an example of a positive $d\times d$ matrix for which the positive stabilizing 
problem has at least $2^d$ local minima. This may be an argument for the high algorithmic complexity of the problem. 
\smallskip 

The paper is organized as follows. In Section~3 we solve the positive 
destabilizing problem and find  the closest non-negative unstable matrix. 
Section~4 deals with the positive stabilizing problem. We present the iteration algorithm for computing the  local minima and prove its convergence with a linear rate. 
We show that if it converges to a positive matrix, then that matrix gives the global minimum. In Section~5 we analyse possible number of local minima and for each $d$, 
give an example of $d\times d$ matrix for which this number is at least $2^d$. 
Then we apply our results for finding closest Hurwitz stable and unstable matrices (Section~6).     

\smallskip

Finally, let us note that in the problem of finding the
closest stable non-negative matrix to a
matrix~$A$, the matrix $A$ itself does not have to be
non-negative. For any real-valued matrix $A$, this problem
can be reduced to the case of non-negative~$A$ by 
considering the matrix $A_{+} =\max \, \{A, 0\}$ (the entrywise
maximum). It is shown easily that  closest stable  matrix to
the matrices $A$ and $A_{+}$ are the same. 
Therefore, in what follows we assume everywhere the initial matrix
$A$ is non-negative.
\smallskip

\begin{center}
\textbf{2. Framework}
\end{center}

Let us first  introduce necessary notation. 
The Frobenius norm of any rectangular matrix $X$ is  
$\|X\| \, = \, \sqrt{{\rm tr}\, X^T X}\, = \, \sqrt{\sum_{i, j}|x_{i, j}|^2}$. This is the standard Euclidean norm in the space of matrices regarded as vectors 
 with the scalar product $\bigl(X, Y \bigr)\, = \, {\rm tr}\, X^T Y$. 
We write $X\perp Y$ for the case $\bigl(X, Y \bigr)\, = \, 0$. We keep the 
notation $\|X\|$ for the Frobenius norm and by $\|X\|_2 = \sqrt{\rho(X^TX)}$
denote the Euclidean operator norm of the matrix~$X$. 

For arbitrary rectangular matrices, $A$ and $B$ for which the product 
$AB$ is well defined we have $\|AB\| \le \|A\|\, \|B\|$. 
In case when $A$ is a co-vector and $B$ is a vector, this becomes an equality. 
If $A$ is a matrix and $B = \bx$ is a vector, this inequality implies 
$$
\|A\bx\|_2 \ = \ \|A\bx\| \ \le \ \|A\|\, \|\bx\|\ = \ \|A\|\, \|\bx\|_2. 
$$
Thus, $\|A\| \, \ge \, \|A\bx\|_2/\|\bx\|_2$ for all $\bx \in \re^d\setminus \{0\}$, 
therefore $\|A\| \ge \|A\|_2$. 

By the Perron-Frobenius theorem, a non-negative matrix possesses 
a non-negative eigenvector corresponding to a non-negative eigenvalue which is bigger than or equal to 
moduli of all other eigenvalues.  This eigenvalue is  called {\em leading}, and the corresponding eigenvector is the {\em leading eigenvector}. 

By the {\em support}
of a non-negative matrix (${\rm supp}\, X$) we mean the set of positions of its positive entries. 

\bigskip 

\begin{center}
\textbf{3. The closest unstable matrix}
\end{center}
\bigskip

For a given non-negative matrix $A$ with $\rho(A) < 1$, we consider the problem 
\begin{equation}\label{eq.prob-m}
\left\{
\begin{array}{l}
\|X \, - \, A\|\ \to \ \min\, \\[2mm]
\rho(X)\ \ge \ 1\, . 
\end{array}
\right.
\end{equation} 
As we show below in Theorem~\ref{th.5}, the solution~$X$ of this problem is always 
a non-negative matrix. Therefore, the search for a closest unstable {\em non-negative} matrix  has the same result. Thus, 
{\em for a non-negative matrix, 
the destabilizing problem is equivalent to the positive destabilising problem.} 
We begin with some auxiliary notation and facts. 

\bigskip 

\begin{center}
\textbf{3.1. Auxiliary facts}
\end{center}
\bigskip

\smallskip 
 
For a given $d\times d$ matrix $A$ such that $\rho(A)\ne 1$ we 
denote 
\begin{equation}\label{eq.MN}
M \ = \ (I-A^T)(I-A)\ ; \qquad N \ = \ (I-A)(I-A^T)
\end{equation}
Since  $M$ and $N$ are both symmetric and positive definite, all their eigenvalues are non-negative.  
Take an arbitrary eigenvector $\bv$ of $M$ associated to some eigenvalue $\mu \ne 0$ and normalize it as $\|\bv\| = 1$. 

Choosing $\, r \, \in \, \bigl\{\sqrt{\mu}\, , \, - \sqrt{\mu}\bigr\}$, 
we denote 
\begin{equation}\label{eq.X}
\bu = \frac{1}{r}\, (I-A)\bv \quad \mbox{and} \quad X \ = \ A \ + \ r\, \bu \, \bv^T\, . 
\end{equation}

\begin{lemma}\label{l.5}
The vector $\bu$ defined in~\eqref{eq.X} satisfies $N\bu = \mu \bu$  and $\|\bu\| = 1$. 
Moreover, $\bv$ and $\bu$ are the right and the left eigenvectors respectively of $X$ 
associated to the eigenvalue one. 
\end{lemma}
{\tt Proof}. Combining the equality $(I-A^T)(I-A)\bv = r^2\bv$ with the 
definition of $\bu$ we get 
 \begin{equation}\label{eq.system}
\left\{
\begin{array}{lll}
(I-A)\ \  \bv & = & r\, \bu \\[2mm]
(I-A)^T\bu & = & r\, \bv\, . 
\end{array}
\right. 
\end{equation}
Since $\|(I-A)\bv\|^2 \, = \, \bigl( (I-A)(I-A^T)\bv \, , \, \bv\bigr)\, = \, 
r^2(\bv, \bv)\, = \, r^2$, taking into account the first equation in~(\ref{eq.system})
we see that $\|r\, \bu\|^2\ = \, r^2$ and hence $\|\bu\| = 1$. 
Furthermore, 
$$
N\bu \ = \ (I-A)(I-A^T)\bu\ = \  (I-A) r\bv \ = \ r^2\bv\, .
$$ 
Finally, $X\bv  \, = \,    A\bv \, + \, r\bu (\bv,  \bv)$. Since $(\bv, \bv) = 1$ and 
$A\bv + r\bu \, = \, \bv$, which  follows from the first equation of~(\ref{eq.system}), we obtain $X\bv = \bv$. In the same way one shows that $\bu^TX = \bu^T$. 

{\hfill $\Box$}
\medskip 

\bigskip 

\begin{center}
\textbf{3.2. A formula for the closest unstable matrix}
\end{center}
\bigskip 

The following theorem provides an explicit solution to the non-negative destabilization problem. 

\begin{theorem}\label{th.5}
Let  $A $ be an arbitrary non-negative matrix  such that $\rho(A) < 1$
and let $M  =  (I-A^T)(I-A)$. 
Let~$\mu$ be the smallest eigenvalue of~$M$. 
Then~$M$ possesses a non-negative eigenvector~$\bv$ associated to $\mu$; 
moreover, for  $r = \sqrt{\mu}$, both the vector $\bu$ and the matrix~$X$ 
defined by~\eqref{eq.X} are non-negative. 
Finally the matrix  $X$ is the closest unstable matrix to $A$ and $\|X-A\| = r$. 
\end{theorem}

Thus, to find the closest unstable matrix to a non-negative matrix~$A$ one needs to 
take the smallest singular value $r$ of the matrix $I-A$ and take the corresponding 
normalized singular vector~$\bv$. There may be a subspace of such vectors, in the case 
when $r$ is multiple, but it always contains a non-negative singular vector~$\bv\ne 0$, 
as it is guaranteed by Theorem~\ref{th.5}. 
Take this vector and normalize it as $\|\bv\| = 1$. 
Then the solution $X$ is readily available by formula~(\ref{eq.X}). 
By Theorem~\ref{th.5}, the matrix~$X$ is non-negative and is the closest unstable matrix 
to $A$ among all matrices, not only non-negative ones. 
\smallskip 

{\tt Proof of Theorem~\ref{th.5}}. Observe that 
\[
(I-A)^{-1} \ = \ \sum_{k=0}^{\infty}A^k \ \ge 0 \ ; 
\qquad (I-A^T)^{-1} \ = \ \sum_{k=0}^{\infty}(A^T)^k \ \ge 0 \ 
\]
(both those series converge since $\rho(A^T) = \rho(A) < 1$). 
Therefore, the matrix $M^{-1}$ is non-negative as well. 
Consequently, its biggest by modulo eigenvalue is 
non-negative and is realized with a non-negative eigenvector~$\bv$. 
The reciprocal to this eigenvalue is the smallest by modulo non-negative 
eigenvalue of $M$. Denote this eigenvalue by $\mu$ and take $r = \sqrt{\mu}$. 
The second equation of the system~(\ref{eq.system}) yields 
$\, \bu \,  = \, r(I-A^T)^{-1}\bv$, hence $\bu \ge 0$, because 
$(I-A)^{-1} \ge 0$. 
Therefore, $X \, = \, A \, + \, r\, \bu\, \bv^T \, \ge \, 0$.
Moreover, Lemma~\ref{l.5} implies that $\|\bu\| = 1$.  
Hence 
\[
\|X - A\| \ = \ r \, \|\bu\, \bv^T\|\ = \ r \, \|\bu\|\, \|\bv^T\| \ = \ r\, . 
\]   
On the other hand, for every matrix $Y$ with $\rho(Y) = 1$, we have 
$\|Y - A\|\ge r$, which proves the optimality of $X$. 
To show this we first assume that  $Y\ge 0$.  In this case $Y$ has a leading eigenvector $\bz$, 
for which $Y\bz = \bz$ and $\|\bz\| = 1$. Then 
$$
\|Y - A\| \ = \ \|Y - A\|\, \|\bz\|\ \ge \  \|(Y - A)\bz\|\ = \ 
\|\bz - A\bz\| \ = \ \|(I - A)\bz\|\ = \ \|(I - A)\bz\|_2\, . 
$$
However, since all singular values of the matrix $I-A$ are bigger than 
or equal to $r$, it follows that $\|(I - A)\bz\|_2 \, \ge \, r\, \|\bz\|_2 = r$. 
Thus, $\|Y - A\|\ge r$, which proves that the matrix $X$ is the closest stable 
matrix for $A$ among non-negative matrices. 

Take now an arbitrary matrix $\Delta$ such that $\|\Delta\| = r$
and show that $\rho(A+\Delta) \le 1$. This will prove the optimality 
of $X$ among all matrices. The matrix $|\Delta|$ composed by the moduli of the entries of $\Delta$ has the same norm $r$. On the other hand, 
$\|(A+|\Delta|)^k\| \, \ge \, \|(A+\Delta)^k\|$ for every $k$, which    
in view of Gelfand's formula for the spectral radius implies that  
$\rho(A+|\Delta|)\ge \rho(A+\Delta)$. On the other hand, since 
$X$ is optimal among non-negative matrices, we see that 
$\rho(A+|\Delta|) \le \rho(X) = 1$, and therefore $\rho(A+\Delta)\le 1$, which 
completes the proof.  
{\hfill $\Box$}
\medskip

\begin{remark}\label{r.38}
{\em Note that since the difference $X-A$ has rank one, $X$ is also 
the closest stable matrix to $A$ in the spectral (Euclidean) norm.}
\end{remark}

\begin{remark}\label{r.40}
{\em The simplicity of  solution of the non-negative 
destabilization problem is explained by the fact this problem is 
actually unconstrained. Indeed, By Theorem~\ref{th.5}, if a matrix $A$ is non-negative, 
then its closest unstable matrix is also non-negative. Hence, 
the non-negativity constraints turn out to be redundant here. On the other hand, 
we can exploit all advantages of the non-negativity of the solution 
provided by the Perron-Frobenus theory. In contrast, in the stabilization 
problem the non-negativity constraints are significant, which makes that problem 
much more difficult. 
}
\end{remark}
\bigskip 

\begin{center}
\textbf{3.3.  Illustrative example}
\end{center}
\bigskip 

Consider the matrix
\[
A = \left( \begin{array}{ccc}
0.4 & 0.4 & 0.1 \\ 
0.5 & 0.3 & 0.3 \\ 
0.1 & 0.1 & 0.5
\end{array}
\right), \qquad \mbox{with} \quad \rho(A)=0.8960.
\]
The minimal eigenvalue of the matrix $M$ is $0.0102$, which gives $r=0.1009$; the
computation of the vectors $\bv$ and $\bu$ gives (to a five digit precision):
\[
\bu = \left( \begin{array}{r} 0.6484  \\  0.5452  \\  0.5314 \end{array} \right)  
\quad \mbox{and} \quad
\bv = \left( \begin{array}{r} 0.6275  \\  0.6852  \\  0.3698 \end{array} \right).
\]
This yields
\[
X = A + r\,\bu\,\bv^T = 
\left( \begin{array}{ccc}
0.4410  &  0.4448  &  0.1242 \\
0.5345  &  0.3377  &  0.3203 \\
0.1336  &  0.1367  &  0.5198
\end{array} \right).
\]
Theorem~\ref{th.5} yields that~$X$ is the closest unstable matrix, i.e., 
provides a global minimum to the destabilizing problem.  Note that applying the general purpose algorithm in \cite{GM} to compute the stability radii,
one gets, as expected, the same matrix~$X$.
\bigskip

\begin{center}
\textbf{4. The closest stable matrix}
\end{center}
\bigskip

For a given 
non-negative matrix $A$ with $\rho(A) > 1$, we consider the problem 
\begin{equation}\label{eq.prob}
\left\{
\begin{array}{l}
\|X \, - \, A\|\ \to \ \min\, \\[2mm]
\rho(X)\ \le \ 1\, , \quad X \, \ge \, 0. 
\end{array}
\right.
\end{equation}  
We indicate by ${\rm locmin}$ the set of local minima for problem \eqref{eq.prob}.
Simple examples show that 
the constraint $X\ge 0$ is significant here. The reason is that even if a matrix $A$ is positive, then its closest stable matrix (in the space of all matrices) may have some negative elements as the following example demonstrates: 
\begin{ex}\label{ex.5}
{\em For  the matrix 
$A= \left( \begin{array}{cc} 2 & 2 \\ 0 & 0 \end{array} \right)$,  
the closest non-negative stable matrix is 
$X= \left( \begin{array}{cc} 1 & 2 \\ 0 & 0 \end{array} \right)$. 
Indeed, for any other non-negative matrix $Y$, we have 
$\rho(Y) \, \ge \, \left( \begin{array}{cc} y_{11} & y_{12} \\ 0 & 0 \end{array} \right) \, = \, y_{11}$. Hence, if $\rho(Y) \le 1$, then $y_{11} \le 1$, and consequently $\|Y - A\| \ge |y_{11} - 2| \ge 1$. Thus, 
$\|Y - A\| \, \|X - A\|$, and $X$ is the closest stable} 
non-negative {\em  matrix to~$A$. 
On the other hand,  
there exists a closer stable matrix, which is not non-negative:  
$Y' = \left( \begin{array}{cc} 2 & 2 \\ -1/2 & 0 \end{array} \right)$, 
for which $\|Y' - A\| = \frac12$.
Hence, for the matrix~$A$, the closest stable  matrix 
is not non-negative. } 
\end{ex}

We see that the stabilization of the matrix can set some of its entries to zero, and this set of zeros can influence the spectral properties of~$X$. That is why, in problem~(\ref{eq.prob}) 
the combinatorics of the matrix~$X$ plays a role. It can be primitive, imprimitive, irreducible, reducible, etc. All this properties have to be considered. This explains the algorithmic complexity of the problem. As we will see in Section~5, problem~(\ref{eq.prob}) may have exponentially many local minima, all with different combinatorics. Explicit solutions can still be obtained, but under special assumptions (Subsection~4.1). In general, we may hope only for algorithmic solutions of finding local minima. This problem requires some preparation;  
we begin with some simple observations.  
\smallskip 

\begin{lemma}\label{l.8}
Suppose $A$ is a reducible matrix, i.e., there exists a permutation matrix 
$\Pi$ which factorizes $A$ to block upper triangular form, 
\[
\Pi A \Pi^T = \wA =  
\left( \begin{array}{ccccc}
\wA_{1,1} & \wA_{1,2} & \ldots & \ldots & \wA_{1,m} \\
{\bf 0} & \wA_{2,2} & \ldots & \ldots & \wA_{2,m} \\
\vdots & & \ddots & & \vdots \\
{\bf 0} & \ldots & {\bf 0} & \wA_{m-1,m-1} & \wA_{m-1,m} \\
{\bf 0} & {\bf 0} & \ldots & {\bf 0} & \wA_{m,m}
\end{array}
\right) 
\]

Then  the closest stable non-negative matrix~$\widetilde X$ to~$\widetilde A$ is given as follows:  
\begin{equation}
X = \Pi^T \wX \Pi \quad \mbox{with} \quad
\wX = \left( \begin{array}{ccccc}
\wX^{(1)} & \wA_{1,2} & \ldots & \ldots & \wA_{1,m} \\
{\bf 0} & \wX^{(2)} & \ldots & \ldots & \wA_{2,m} \\
\vdots & & \ddots & & \vdots \\
{\bf 0} & \ldots & {\bf 0} & \wX^{(m-1)} & \wA_{m-1,m} \\
{\bf 0} & {\bf 0} & \ldots & {\bf 0} & \wX^{(m)}
\end{array}
\right) 
\label{eq:X}
\end{equation}
where $\wX^{(i)}$ is the closest stable non-negative matrix to $\wA_{i,i}$ 
(for $i=1,\ldots,m$).
In particular if $\rho(\wA_{i,i}) \le 1$, then 
$\wX^{(i)}=\wA_{i,i}$.
\end{lemma}
{\tt Proof.} 
Since $\wA$ is similar to $A$ and $\Pi$ defines an isometry in the Frobenius norm,
considering the problem of finding the closest stable matrix $X$ 
to $A$ is equivalent to that of finding the closest stable matrix 
$\wX$ to $\wA$.

Let $\wX$ be the matrix constructed in \eqref{eq:X}. Since each block 
$\wX^{(i)}$ is stable, we have $\rho(\wX) = \max_i \, \rho(\wX^{(i)}) = 1$, hence $\wX$
is stable. Any change of some element off the diagonal blocks would increase the distance 
to~$\wA$. On the other hand, this change would not reduce the spectral radius of $\wX$. 
Finally, any change of a diagonal block which keeps this block stable 
would not reduce the distance to~$\wA$. Hence, $\wX$ gives the global minimum. 
{\hfill $\Box$}

\medskip 

Lemma~\ref{l.8} reduces problem~(\ref{eq.prob}) to several similar problems of smaller dimensions whenever $A$ is reducible. Hence, we do not consider this case 
any more. 
{\em In the sequel  we assume that $A \ge 0, \, \rho(A) < 1$, and that 
$A$ is irreducible.} 

\begin{lemma}\label{l.10}
If $X$ is a local minimum for~\eqref{eq.prob}, then $X \le A$ and 
$\rho(X) = 1$.  
\end{lemma}
{\tt Proof.} If $X$ has a component $x_{ij}$ bigger than $a_{ij}$, then we slightly reduce~$x_{ij}$. 
The distance $\|X - A\|$ decreases, while the spectral radius $\rho(X)$, as a monotone 
function on the set of non-negative matrices, does not increase. Hence the condition $\rho(X)\le 1$
remains true, which contradicts to the local optimality of~$X$. If $\rho(A) > 1$, then $X \ne A$, 
and hence $x_{ij} < a_{ij}$ at least for one component. If $\rho(X) < 1$, then we can slightly increase 
$x_{ij}$ so that the condition $\rho(X)\le 1$
remains true. 

{\hfill $\Box$}
\medskip 

Lemma~\ref{l.10} ensures that we can restrict our search to 
matrices that are entrywise smaller or equal to~$A$ and have spectral radius one. 
So,  in the sequel we assume that $X \le A$. Thus, we look for solutions of the problem~(\ref{eq.prob})
on the set of matrices $X$ such that $0 \, \le \, X \, \le \, A$ and $\rho(X) = 1$. 
\smallskip 

Thus, the leading eigenvalue of $X$ is equal to one.
If this  eigenvalue is simple, in particular, if the matrix $X$
is irreducible, then the spectral radius is differentiable at 
the point~$X$, and problem~(\ref{eq.prob}) is smooth. 
Then we can apply the Lagrange theorem and derive the following condition for the local minimum: 
\begin{prop}\label{p.10}
Suppose a matrix $X$ is a local minimum for~\eqref{eq.prob} and its 
leading eigenvalue~$1$ is simple. Suppose $\bu$ and $\bv$ are respectively left and right leading eigenvectors of~$X$,  associated to the eigenvalue $1$;   then there exists  a number $r > 0$
and a matrix $\Lambda \ge 0$
 such that $\Lambda \perp X$ and 
\begin{equation}\label{eq.kkt}
A \quad = \quad X \ + \ r\, \bu \, \bv^T \ - \ \Lambda\, . 
\end{equation}   
\end{prop}
\begin{remark}\label{r.18}{\em 
Under the assumptions of the proposition, we have 
$(\bv , \bu) > 0$. Indeed, if $(\bv , \bu) = 0$, then 
the supports of the vectors $\bv$ and $\bu$ are disjoint. 
After a permutation of the basis  the matrix $X$ obtains the block diagonal form 
$X \, = \, \left( \begin{array}{cc} X^{(1, 1)} & 0 \\ 0 & X^{(2, 2)} \end{array} \right)$   with the blocks corresponding to the supports of $\bv$ and $\bu$ respectively and with $\rho(X^{(1, 1)}) = \rho(X^{(2, 2)})= 1$. Hence $1$ is not a 
simple eigenvalue. 
}
\end{remark}

{\tt Proof.} The derivative of the function $\|X-A\|^2$ at $X$ is equal to 
$2\,(X - A)$. Since the leading eigenvalue~$\lambda_{\max}$ of~$X$ is simple, 
the function $\lambda_{\max}(X)$ is differentiable at~$X$ and  the 
 gradient 
is equal to $\bu\,  \bv^T$ (see~\cite{K}). Now applying the Lagrange theorem to 
the problem~(\ref{eq.prob}) we see that there are non-negative multipliers $\alpha_0, \alpha_1$ and 
$\alpha_{ij},  i ,j = 1, \ldots , d$, such that 
$$
2\, \alpha_0\, (X \, - \, A)\ + \ \alpha_1 \, \bu \, \bv^T \ - \ \Lambda\ = \ 0, 
$$ 
where $\Lambda \, = \, (\alpha_{ij})_{i, j}$ is the matrix of multipliers corresponding to 
the constraints $x_{i j} \ge 0$. The complementary slackness conditions 
give $\alpha_{ij}x_{ij} = 0$ for all $i, j$, and hence $\Lambda \perp X$. 

If $\alpha_0 = 0$, 
then $ \bu \, \bv^T \, \perp \, X$.  In this case 
$$
0 \ = \ \bigl( \bu \, \bv^T \, , \, X \bigr) \ = \ 
{\rm tr} \, \bigl( (\bu \, \bv^T)^T X \bigr) \ = \ 
{\rm tr} \, \bigl( \bv \, \bu^T X \bigr)\, . 
$$
Since $\bu^T X = \bu^T$, we see that 
${\rm tr} \, \bigl( \bv \, \bu^T X \bigr) \, = \,  {\rm tr} \, \bigl( \bv \, \bu^T\bigr) \ = \ (\bv , \bu). $ 
Thus, if $\alpha_0 = 0$, then  $\bv \perp  \bu$. This means that 
the leading eigenvalue of $X$ is multiple, which contradicts the assumption.  
 Thus,  $\alpha_0 > 0$, and we can set $\alpha_0 = \frac12$ and arrive at~(\ref{eq.kkt}).  
   
{\hfill $\Box$}
\medskip

Equation~(\ref{eq.kkt}) is not simple to solve, because it involves an unknown matrix~$X$
together with its left and right eigenvectors. Nevertheless, if it possesses a positive solution~$X$, then it can be found explicitly (Subsection~4.1). 
In general, $X$ can have zero entries, and therefore 
additional unknowns occur: each zero entry of~$X$ generates the corresponding 
 unknown element of the matrix~$\Lambda$. 
In Subsection~4.2 we present an algorithm for the numerical solution of problem~(\ref{eq.kkt}). However, in some cases this solution is not able to identify a point of local minimum in a unique way. This happens when the matrix $X$
is non-primitive. Moreover, if the eigenvalue~$1$ is multiple for~$X$, 
then Proposition~\ref{p.10} may not hold at all. In this case  
 the matrix $X$ must be reducible (see, for instance,~\cite[chapter 13, \S 2, theorem 2]{G}), and this case is considered 
 in the end of this section. 
 
We see that the sparsity pattern of the matrix~$X$, i.e., the location of zero components, 
is crucial in the solution of~(\ref{eq.kkt}), because it defines the set of extra variables 
in the matrix~$\Lambda$. That is why the solution involves the combinatorics of the matrix~$X$. 

Let us recall that a matrix $X \ge 0$ is called {\em primitive} if some of its power is strictly positive. If~$X$ is non-primitive, but irreducible, it is called 
{\em imprimitive}. We have to analyse conditions for the local minimum 
in the three separate cases: 1) $X$ is primitive; 2) $X$ is imprimitive (a quite unusual case in our experiments);
3) $X$ is reducible. 
\smallskip 

\noindent \textbf{Case 1. $\mathbf{X}$ is primitive}. We call a primitive  matrix $X$ satisfying conditions of Proposition~\ref{p.10} a {\em stationary point} of problem~\eqref{eq.prob}. 
\smallskip

\noindent \textbf{Case 2. $\mathbf{X}$ is imprimitive}.
In this case one more necessary condition to local minimality appears. By the Perron-Frobenius theorem, 
for an imprimitive matrix~$X$, there is a disjoint partition of the set $\Omega = \{1, \ldots , d\}$ into  
$r\ge 2$ nonempty sets $\Omega_1, \ldots , \Omega_r$ such that the matrix $X_k$ defines 
a cyclic permutation  of those sets: $\Omega_1 \to \Omega_2\to \cdots \to \Omega_r \to \Omega_1$. This means that if ${\rm supp} \, \ba \, \subset \, \Omega_i$, then 
${\rm supp} (X\ba) \, \subset \, \Omega_{i+1}, \, 
i = 1, \ldots , r$ (we set $\Omega_{r+1} = \Omega_1$). 
After renumbering the basis vectors, $X$ gets the form 
of cyclic permutation of  primitive blocks $X^{[1]}, \ldots , X^{[r]}$. 
\begin{equation}\label{eq.cyclic}
X \quad = \quad \left(
\begin{array}{cccccc}
0 & 0 &  0 & \ldots &  0 & X^{[r]}\\
X^{[1]} & 0 &  0 & \ldots & 0 & 0\\
0 & X^{[2]} &  0 & \ldots & 0 & 0\\
\vdots & {} &  {} & \vdots  & {} & \vdots \\
0 & 0 & 0 & \ldots & X^{[r-1]} & 0
\end{array}
\right)\ .
\end{equation}
For a positive {\em vector of weights} $(s_1, \ldots , s_r)$, denote by $X[s_1, \ldots , s_r]$
the same matrix with blocks $s_1X_k^{[1]}, \ldots , s_rX_k^{[r]}$. 
The spectral radius of this matrix is equal to 
$s_1\cdots s_r\rho(X)$. Then we optimize the weights~$s_1, \ldots , s_n$ by solving the problem
 \begin{equation}\label{eq.prim}
\left\{
\begin{array}{l}
\|X[s_1, \ldots , s_r] \, - \, A\|\ \to \ \min\\
s_1\cdots s_r\ = \ 1 
\end{array}
\right.
\end{equation} 
This problem always admits a unique point of minimum, the details of the solution are outlined  to Sections~4.3. Now we can extend the notion of stationary point to all 
irreducible matrices. 

\begin{defi}\label{d.12}
An irreducible matrix is stationary for problem~\eqref{eq.prob}
if it satisfies the equation~\eqref{eq.kkt} and, if it is imprimitive, 
has the optimal weights of the blocks $s_1, \ldots , s_m$ obtained by 
solving problem~\eqref{eq.prim}.  
\end{defi}
The case $r=1$ corresponds to a primitive matrix, the case $r \ge 2$ 
does to imprimitive one. 
\smallskip 

\noindent \textbf{Case 3. $\mathbf{X}$ is  reducible}. 
In this case the matrix $X$ admits a unique, up to 
a permutation of the basis vectors, {\em Frobenius form}. This means that
there exists a reordering  of the basis vectors, after which $X$
gets the following block upper-triangular form:  
\begin{equation}\label{eq.frob}
X\quad = \quad
\left(
\begin{array}{cccc}
 X^{(1)} & * & \cdots & *\\
 \bzero       & X^{(2)} &   *   & \vdots \\
\vdots   &         & \ddots & * \\
\bzero & \cdots  & \bzero & X^{(m)}
\end{array}
\right).
\end{equation}
where all matrices $X^{(i)}$ in the diagonal blocks are irreducible. 
For an irreducible matrix~$X$, we have $m=1$, otherwise $m\ge 2$.
Now we can define the notion of stationary point for a general non-negative matrix~$X$.   
\begin{defi}\label{d.16}
A matrix~$X$ is said stationary for problem~\eqref{eq.prob} if its Frobenius 
form~\eqref{eq.frob} is such that above the diagonal blocks we have $X=A$, and 
for the diagonal blocks $X^{(i)}$ we have:  
\begin{itemize}
\item[(1) ] if $\rho(A^{(i)})\le 1$, then $X^{(i)} =  A^{(i)}$; 

\item[(2) ] if $\rho(A^{(i)})> 1$, then $\rho(X^{(i)}) = 1$ and 
$X^{(i)}$ is a stationary point 
(by Definition~\ref{d.12}) of the problem  
$\|X^{(i)} -  A^{(i)}\| \to \min, \ \rho(X^{(i)})\le 1$. 
\end {itemize}
\end{defi}
If  $m=1$, then this definition is reduced to the cases 1) and 2). 
\begin{theorem}\label{th.12}
If $X \in {\rm locmin}$ for problem~\eqref{eq.prob}, then $X$ is stationary. 
\end{theorem}
{\tt Proof.} Consider the Frobenius form of the matrix $X$. 
If some element $x_{ij}$ over the diagonal block is not equal 
to $a_{ij}$, then $x_{ij} < a_{ij}$ (Lemma~\ref{l.10}). If we 
slightly increase this element, the distance $\|X-A\|$ is reduced while 
the spectral radius of $X$ does not change, which contradicts the local 
optimality of~$X$.  Hence, $X=A$ beyond the diagonal blocks. 
Consider now any diagonal block~$X^{(i)}$. If $\rho(A^{(i)})\le 1$, 
but $X^{(i)}\ne A^{(i)}$, then we consider a matrix 
$X^{(i)}(t) =  (1-t)\, X^{(i)} \, + \, t\, A^{(i)}$. Since $X^{(i)} \le A^{(i)}$, 
all entries of the matrix $X^{(i)}(t)$ increase in~$t$, and therefore so 
does the spectral radius $\rho(X^{(i)}(t))$. Hence, for every $t\in  (0, 1)$, 
the spectral radius of the matrix $X^{(i)}(t)$ does not exceed one, while 
this matrix is closer to~$A^{(i)}$ than~$X^{(i)}$. This again contradicts the 
local optimality of~$X$. Finally, if $\rho(A^{(i)}) > 1$, then 
$X^{(i)}$ must be the local minimum of the problem  
$\|X^{(i)} -  A^{(i)}\| \to \min, \ \rho(X^{(i)})\le 1$. 
Since $X^{(i)}$ is irreducible, by Proposition~\ref{p.10} it satisfies 
equation~(\ref{eq.kkt}) and, if it is imprimitive, 
has the optimal weights of the blocks $s_1, \ldots , s_m$ obtained by 
solving problem~(\ref{eq.prim}). 
Hence, $X^{(i)}$ is stationary.

{\hfill $\Box$}
\medskip

Thus, if a local minimum is attained at a matrix~$X$, then $X$ is stationary:  if $X$ is primitive, then it satisfies equation~(\ref{eq.kkt}), if it is imprimitive, 
then it also satisfies an additional optimality  condition~(\ref{eq.prim}), if it is 
reducible, then its Frobenius form consists of stationary blocks (primitive or imprimitive) and satisfies the requirements stated in Definition~\ref{d.16}. In Subsection~4.2 we will see how to construct the stationary matrices and to find local minima of problem~(\ref{eq.prob}) algorithmically. 

We begin with a special case when the 
stationary matrix is strictly positive. In this case it satisfies equation~(\ref{eq.kkt}) with $\Lambda = 0$. It turns out that under some extra assumptions it provides the global minimum to problem~(\ref{eq.prob}).

\bigskip 
 

\bigskip 

\begin{center}
\textbf{4.1. A positive local minimum is a  global minimum}
\end{center}

\bigskip 

We show here that if a strictly positive matrix gives a local minimum to problem~(\ref{eq.prob}), then it gives its global minimum. 
Moreover, this matrix can be explicitly found.  

Consider the solution of the {\em destabilization} problem in Section~3. 
Can we apply the same reasoning to the stabilization problem, assuming  $\rho(A) > 1$?  
We define the matrices $M$ and $N$ by the same formula~(\ref{eq.MN}), then 
define $\bv$ as an eigenvector of $M$ corresponding to its smallest eigenvalue, then define $\bu$ and $X$ by the same formula~(\ref{eq.X}) with $r=-\sqrt{\mu}$ instead of 
$\sqrt{\mu}$ (this is the only difference!).  Then we  repeat the proof of 
Theorem~\ref{th.5} to establish that $X$ is the closest stable matrix. However, 
here we cannot show the positivity of the matrix  $M^{-1}$  (this already may not be true) and hence, the positivity of $X$.  

Thus, in general,  Theorem~\ref{th.5} cannot be extended to the stabilization problem. 
Nevertheless, if the obtained vectors $\bv, \bu$ and the matrix $X$ are non-negative, then 
$X$ is the global minimum for problem~(\ref{eq.prob}), and the corresponding proof is literally the same as the proof of Theorem~\ref{th.5}. We formulate it in the following  

\begin{theorem}\label{th.7}
Assume the matrix~$M$ (see \eqref{eq.MN}) possesses a non-negative eigenvector $\bv$, 
corresponding to its smallest eigenvalue~$\mu$. 
Assume also that for  $r = - \sqrt{\mu}$, 
the vector $\bu$ and the matrix~$X$ defined by~\eqref{eq.X} are both non-negative. 
Then  $X$ is the closest stable non-negative matrix for $A$ and $\|X-A\| = r$.  
\end{theorem}

\medskip 

Now we formulate the main result of this subsection. 

\begin{theorem}\label{th.8}
If a matrix $X > 0$ provides a local minimum for problem~\eqref{eq.prob}, 
then it provides a global minimum. Moreover, in this case all assumptions of 
Theorem~\ref{th.7} are fulfilled and $X$ coincides with the corresponding matrix 
from that theorem.  
\end{theorem}
{\tt Proof}. Each point of local minimum has the form~(\ref{eq.kkt}), where 
in the case $X>0$, we have $\Lambda = 0$. Thus, $X \, = \, A - r\, \bu \, \bv^T$. 
Multiplying this equality by $\bv$ from the right and taking into account that 
$\bv^T \bv \, = \, (\bv, \bv ) = 1$ and that $X\bv = \bv$, we obtain 
$(I-A)\bv \, = \, - r\, \bu$ and hence $r^2 = \|(I-A)\bv\| \, = (\bv, M\bv)$. 
On the other hand, $\|X-A\| = r$. Thus, $\|X-A\|^2 \, = \, (\bv, M\bv)$. 
Substituting $\bu = -\frac{1}{r}(I-A)\bv$ in the formula for~$X$ we get 
\begin{equation}\label{eq.xnew}
 X \quad = \quad A \ - \ (A-I)\, \bv\,  \bv^T\, .
\end{equation}
If $\bv$ is not an eigenvector of $M$ corresponding to its smallest eigenvalue, 
then there is a  vector $\tilde \bv$ close to $\bv$ such that $\|\tilde \bv\| = 1$
and $(\tilde \bv, M\tilde \bv)\, < \, (\bv, M\bv)$. Define the matrix 
$$
 \tilde X \quad = \quad A \ - \ (A-I)\, \tilde \bv\,  \tilde \bv^T\, .
$$
If $\tilde \bv$ is close enough to $\bv$, then $\tilde X >0$. Moreover, 
 $\tilde X \tilde \bv = \tilde \bv$,  hence the spectral radius of~$\tilde X$
 is one. Finally, 
 $\|\tilde X - A\|^2 \, = \, \|(A-I)\tilde \bv\|^2 \, = \, 
(\tilde \bv, M\tilde \bv)$, which is smaller than 
$(\bv, M\bv) \, = \, \|X-A\|^2$. Thus, $\tilde X$ is closer to $A$ than $X$, hence 
$X \notin {\rm locmin}$. The contradiction proves that the leading eigenvector~$\bv$ of $X$ is a positive eigenvector of $M$ corresponding to the smallest 
eigenvalue of $M$. By the same argument we show that  the right leading eigenvector~$\bu$ of $X$ is a positive eigenvector of $N$ corresponding to its smallest eigenvalue. Thus, all assumptions of Theorem~\ref{th.7} are fulfilled and $X$ coincides with the corresponding matrix from that theorem.

{\hfill $\Box$}
\medskip 

Thus, problem~\eqref{eq.prob} does not have strictly positive local minima except for those constructed by Theorem~\eqref{th.7}. In particular, we have proved
\begin{cor}\label{c.10}
If the assumptions of Theorem~\eqref{th.7} are not fulfilled, then 
problem~\eqref{eq.prob} does not have strictly positive local minima. 
\end{cor} 

\begin{remark}\label{r.10}
{\em Theorem~\ref{th.8} admits the case when there are infinitely many 
closest stable matrices for~$A$. This happens when the smallest eigenvalue of $M$
is multiple. 

Note that problem~(\ref{eq.prob}) may have positive 
stationary points different from local minima. They have the same form 
$X \, = \, A \, - \, r\, \bu \bv^T$, but with $\mu = r^2$ to be a non-minimal 
eigenvalue of~$M$. On the other hand, problem~(\ref{eq.prob}) may possess 
 positive stationary points corresponding to at most one eigenvalue of~$M$. 
 Indeed, if there are two stationary points 
$X$ and $X'$ corresponding to different eigenvalues $\mu$ and $\mu'$, then 
the corresponding eigenvectors $\bv$ and $\bv'$ are orthogonal to each other. 
On the other hand, they are both non-negative, hence $\bv$ must have at least one zero component. This is impossible, because $\bv$ is an eigenvector of a 
strictly positive matrix~$X$. We collect those observations in the following corollary. 
}
\end{remark}
\begin{cor}\label{c.15}
If problem~(\ref{eq.prob}) possesses  positive stationary points~$X$, they all 
have the form~(\ref{eq.xnew}) with $\bv$ being 
an eigenvector of the matrix~$M$. 
\end{cor}

\begin{remark}\label{r13}
{\em Note that since the difference $X-A$ has rank one in this case $X$ is also 
the closest stable matrix to $A$ in the spectral norm.}
\end{remark}

\begin{ex}\label{ex.33}
{\em Consider the matrix
\[
A = \left( \begin{array}{ccc}
0.6 & 0.4 & 0.1 \\ 
0.5 & 0.5 & 0.3 \\ 
0.1 & 0.1 & 0.7
\end{array}
\right), \qquad \mbox{with} \quad \rho(A)=1.0960.
\]
The minimal eigenvalue of the matrix $M$ is $0.0082$, which gives $r=0.0903$; the
computation of the vectors $\bv$ and $\bu$ gives (to a five digit precision):
\[
\bu = \left( \begin{array}{r} 0.6193  \\  0.4888  \\  0.6144 \end{array} \right)  
\quad \mbox{and} \quad
\bv = \left( \begin{array}{r} 0.6438  \\  0.7166  \\  0.2684 \end{array} \right).
\]
This yields
\[
X = A - r\,\bu\,\bv^T = 
\left( \begin{array}{ccc}
0.5640  &  0.3599  &  0.0850 \\
0.4716  &  0.4684  &  0.2881 \\
0.0643  &  0.0602  &  0.6851 
\end{array} \right).
\]
which is stable and, by Theorem~\ref{th.8}, has minimal distance to $A$.}
\end{ex}

\begin{ex}\label{ex.10}
{\em Fix the dimension~$d\ge 2$ and denote by $E$ the matrix of all ones.
We find the closest stable matrix to the matrix $A = \alpha E$ depending 
on the parameter $\alpha$. 

For $\alpha \le \frac{1}{d}$, the matrix $\alpha E$ is stable, hence $X=\alpha E$ is the global minimum (and a unique local minimum).  
If $\alpha > \frac{1}{d}$, the matrix $M = (A-I)(A^T -I) = (\alpha E - I)^2$
has eigenvalues $(d\alpha - 1)^2, 1 , \ldots , 1$. Hence, if $\alpha \in \
\bigl[\frac{1}{d}, \frac{2}{d}\bigr]$, then $(d\alpha - 1)^2 \le 1$, and therefore
the  eigenvector $\bv = \frac{1}{\sqrt{d}}\be$ is associated  to the smallest 
eigenvalue of~$M$. Hence the matrix $X = \frac{1}{d} E$ provides a 
a global minimum. Indeed, $X = A - r \bu \bv^T$ with 
$\bu = \bv = \frac{1}{\sqrt{d}}\be$ and $r = \bigl(\alpha - \frac{1}{d}\bigr)$. 
It remains to refer to Theorem~\ref{th.7}. However, for each $\alpha > \frac{2}{d}$, 
the vector $\bv$ does not correspond to the smallest eigenvalue of $X$, 
hence, in view of Theorem~\ref{th.7}, the same matrix $X = \frac{1}{d}E$
does not provide even a local minimum. Although this matrix is still 
a stationary point since it is positive and has the form 
$X = A - r \bu \bv^T$. 

Surprisingly, the natural answer~$X = \frac{1}{d}E$ turns out to be wrong for all~$\alpha > \frac{2}{d}$: the matrix $X$ is not even locally closest stable matrix for~$\alpha E$.  
    }
\end{ex}
\begin{ex}\label{ex.20}
{\em Let $\alpha > 0$ and 
$$
A \ = \ \left(
\begin{array}{cc} 
\alpha & \alpha \\
\alpha & \alpha
\end{array}
\right) \ ; \qquad 
X_0 \ = \  
 \left(
\begin{array}{cc} 
\frac12 & \frac12 \\
\frac12 & \frac12
\end{array}
\right) 
$$ 
In view of Example~\ref{ex.10}, if $\alpha \in \bigl[0, \frac12\bigr]$, then 
$X = A$ is the global minimum, if $\alpha \in \bigl[\frac12, 1)$, then 
$X = X_0$ is the global minimum and a unique local minimum. 
For $\alpha  = 1$, the matrix $M = 4I$ has equal eigenvalues, and hence 
$A$ has infinitely many closest stable matrices: every matrix 
\[
X_{t} = A -  \bu_{t} \bv^T_{t}, \qquad \mbox{where} \quad \bu_{t} = 
\bv_{t} = \left(\cos t \ \sin t \right)^T, 
\]
is a global minimum with 
$\|X_t-A\| = 1$. 

Finally, if $\alpha > 1$, then, in view of Corollary~\ref{c.15}, 
$A$ does not have positive local minima, although it has an obvious stationary 
point~$X_0$, which is not locally closest any more. 
Hence, the closest stable matrix has a zero entry. 
Considering two possible cases, when this zero is either off the diagonal
(in this case $X$ is reducible, and hence has ones on the diagonal and $\alpha$
and $0$ off the diagonal) or on the diagonal (this case does not provide minima), 
we conclude that there are two global minima: 
$$
\left(
\begin{array}{cc} 
1 & \alpha \\
0 & 1
\end{array}
\right) \ ; \qquad 
 \left(
\begin{array}{cc} 
1 & 0\\
\alpha & 1
\end{array}
\right) 
$$ 
}
\end{ex}

\bigskip 

\begin{center}
\textbf{4.2. The general relaxation scheme}
\end{center}
\bigskip 

Now we are going to tackle the general case: for a non-negative irreducible 
matrix $A$ such that $\rho(A) > 1$, solve the stabilization problem~(\ref{eq.prob}).  
The idea of the algorithm is the following. We take a matrix $X_0\ge 0$
whose support is not smaller than the support of $A$, and normalize it so that 
$\rho(X_0) = 1$. Then we compute its leading eigenvector $\bv_0$, for which 
$X_0\bv_0 = \bv_0$, and solve the problem 

\begin{equation}\label{eq.prob-1}
\left\{
\begin{array}{l}
\|X \, - \, A\|\ \to \ \min\, \\[2mm]
X\bv_0 \le \bv_0, \ X\ge 0.
\end{array}
\right.
\end{equation} 
Its solution is denoted as $X_1$.

Then we compute the left leading eigenvector $\bu_1$ of $X_1$, for which 
$\bu_1^TX_1 = \bu_1^T$, and solve the problem 
\begin{equation}\label{eq.prob-2}
\left\{
\begin{array}{l}
\|X \, - \, A\|\ \to \ \min\, \\[2mm]
\bu_1^T X \le \bu_1^T, \ X\ge 0.
\end{array}
\right.
\end{equation} 
Its solution is denoted as $X_2$.
Then we loop by alternating between problems of the form \eqref{eq.prob-1} and \eqref{eq.prob-2}, that is
we compute the right leading eigenvector $\bv_2$ of $X_2$ and continue.

To summarize we make a consecutive relaxation of the objective function $\|X-A\|$
every time alternating the right and left leading eigenvectors of the matrix~$X$. 
For even~$k$, we optimize $X_k$ with respect to the fixed leading eigenvector, 
and for odd, we do it with respect to the fixed leading right eigenvector.

We shall prove that if the value $\|X_k-A\|$ is the same for two consecutive iterations, then 
the algorithm halts at the matrix $X_k$. 

Otherwise, the objective function (the distance to $A$) decreases in each iteration. 
The complexity and convergence analysis will be done in the following subsection. 

Before introducing the algorithm, 
we recall some auxiliary facts. If a current matrix $X_{k-1}$
is reducible, then after a suitable renumbering  of the basis vectors it 
gets the form 
 \begin{equation}\label{eq.reduc}
X_{k} \quad = \quad \left(
\begin{array}{cc}
X^{(1, 1)} & X^{(1, 2)} \\
\bzero & X^{(2, 2)}  
\end{array}
\right)\ .
\end{equation}
We denote by $A^{(i, j)}$ the corresponding blocks 
 of the matrix $A$ (after the same renumbering). 
\bigskip 

\begin{center}
\textbf{4.3. The Algorithm}
\end{center}
\bigskip

The Algorithm is based on an inner iteration which implements 
an iterative method to minimize the function $\| X - A \|$ under
the constraints $\rho(X) \le 1$ and $X \ge 0$, and an outer iteration
which takes into account about the reducibility/imprimitivity of
the results provided by the inner iteration and is able to further
refine the construction of a locally optimal solution to the
problem. The outer iteration has a recursive structure, which makes
use of the inner iteration possibly several times, until it halts 
on a stationary point.

\subsection*{The inner iterative optimization Algorithm}

The optimization Algorithm \ref{algo:opt} is a descent method whose flowchart follows.
It constructs a sequence of matrices $\{X_k\}$ which may converge to an irreducible/reducible matrix
and in the first case to a primitive/imprimitive matrix.

\begin{algorithm}[H] \label{algo:opt}
\DontPrintSemicolon
\KwData{$A, X_0$}
\KwResult{$X_{k^+}$, {\bf Reduce}} 
\Begin{
\For{$k=0,\ldots,k_{\max}$}
{
\nl \eIf{$k$ is odd}
{
\nl Compute the right leading eigenvector $\bv_k$ of $X_k$\;
\eIf{$\bv_k > 0$}
{Solve the optimization problem
$\left\{
\begin{array}{l}
\|X \, - \, A\|\ \to \ \min\, \\[2mm]
X\bv_k \le \bv_k, \ X\ge 0.
\end{array}
\right.
$}
{
\nl Set {\bf Reduce}={\bf True}\;
}
}
{
\nl Compute the left leading eigenvector $\bu_k$ of $X_k$\;
\eIf{$\bu_k > 0$}
{Solve the optimization problem
$\left\{
\begin{array}{l}
\|X \, - \, A\|\ \to \ \min\, \\[2mm]
\bu_k^T X \le \bu_k^T, \ X\ge 0.
\end{array}
\right.
$}
{
\nl Set {\bf Reduce}={\bf True}\;
}
}
\nl Set $k^+=k+1$ and $k^-=k-1$\;
\nl Let $X_{k^+}$ be the solution of the optimization problem\;
\nl \If{{\bf Reduce}={\bf True}}
{
\Return
}
\If{$\| X_{k^+}-X_{k} \| < {\rm tol}$ {\bf and} $\| X_{k^-}-X_{k} \| < {\rm tol}$} 
{
\Return
}
}
}
\caption{The iterative optimization Algorithm}
\end{algorithm}


%
%
Basically Algorithm \ref{algo:opt} works as follows.
If $\bv_{k}$ has zeros, then it stops and returns a reducible matrix. 
Otherwise,  if $\bv_{k}>0$, it computes the unique non-negative solution of the problem 
\begin{equation}\label{eq.v}
\left\{
\begin{array}{l}
\|X \, - \, A\|^2\ \to \ \min\\
X\, \bv_{k}\ \le \ \bv_{k} 
\end{array}
\right.
\end{equation} 
and proceeds to next iteration. 

Similarly, if $\bu_{k}>0$, then it computes the unique non-negative solution of the problem
\begin{equation}\label{eq.u}
\left\{
\begin{array}{l}
\|X \, - \, A\|^2\ \to \ \min\\
\bu_{k}^TX\, \ \le \ \bu_{k}^T 
\end{array}
\right.
\end{equation} 
and proceeds to next iteration. 


If $\bv_{k}$ or $\bu_{k}$ has zeros, in fact, the matrix $X_{k}$ is reducible. 

\subsection*{The outer recursive Algorithm}

The outer Algorithm follows.
We take an arbitrary initial matrix~$X_0$ with $\rho(X_0) = 1$ and with the same support of~$A$.
Then we apply the algorithm recursively until a local minimum is found.  
\smallskip 

\begin{algorithm}[H] \label{algo:main}
\DontPrintSemicolon
\KwData{$A, X_0$}
\KwResult{$X$} 
\Begin{
{
\nl Apply Algorithm \ref{algo:opt} with inputs $A$ and $X_0$ and outputs $X$ and {\bf Reduce}\;
\If{{\bf Reduce}={\bf True}}
{Reorder components to get
\[
X  =  \left(
\begin{array}{cc}
X^{(1, 1)} & X^{(1, 2)} \\
\bzero & X^{(2, 2)}  
\end{array}
\right), \qquad
A = \left(
\begin{array}{cc}
A^{(1, 1)} & A^{(1, 2)} \\
A^{(2, 1)} & A^{(2, 2)}  
\end{array}
\right)
\]
\eIf{$\rho(A^{(2, 2)})\, < \,  1$}
{
\nl Set $X^{(2, 2)} \, = \, A^{(2, 2)}$\;
\nl Apply Algorithm \ref{algo:main} with inputs $A^{(1,1)}$ and $X^{(1,1)}$ and output $X^{(1,1)}$
}
{Apply Algorithm \ref{algo:main} with inputs $A^{(1,1)}$ and $X^{(1,1)}$ and output $X^{(1,1)}$\;
 Apply Algorithm \ref{algo:main} with inputs $A^{(2,2)}$ and $X^{(2,2)}$ and output $X^{(2,2)}$
}
}
}
\If{$X$ is imprimitive}
{
find the sets~$\Omega_1, \ldots,\Omega_r$ and compute the optimal weights $s_1, \ldots , s_r$ 
by solving the problem~\eqref{eq.prim}\; 
\nl Set $X = X[s_1, \ldots , s_r]$.
}
\Return
}
\caption{The main recursive Algorithm}
\end{algorithm}

If the matrix $X$ is reducible, it gets the block upper triangular form~(\ref{eq.reduc}), 
where $X^{(1, 1)}$ is an $m\times m$-matrix and $X^{(2, 2)}$ is a $(d-m)\times (d-m)$-matrix 
respectively. 

Consider the case when the right leading eigenvector $\bv$ has zeros, the case of $\bu$ is  similar. 
After renumbering it can be assumed that the first $m$ entries of  $\bv$ are positive and the 
other $d-m$ are zeros. 

Note that $\max \, \bigl\{ \rho(X^{(1, 1)})\, ,  \, \rho(X^{(2, 2)})\, \bigr\}\,  = \, 
\rho(X)\, = \, 1$. Moreover, $\rho(X^{(1, 1)}) \, \ge \, \rho(X^{(2, 2)})$, otherwise, 
the leading eigenvector $\bv$ cannot have zeros in the last $d-m$ positions (in that case
the corresponding eigenvalue would be $\lambda_{\max}(X^{(1, 1)})\, < \, 1$). 
Thus, 
\begin{equation}\label{eq.lambdas}
\rho(X^{(1, 1)}) \ = \ 1\ , \qquad \rho(X^{(2, 2)}) \ \le \ 1\,  .
\end{equation}
 %


Then we first set $X^{(1, 2)} \, = \, A^{(1, 2)}$.  This reduces the distance $\|X - A\|$
and does not change the spectral radius of $X$. 
Second the Algorithm proceeds following one of the following two cases.

{\tt Case Red-1.} $\ \rho(A^{(2, 2)})\, < \,  1$. 
We set $X^{(2, 2)} \, = \, A^{(2, 2)}$. This reduces the distance $\|X - A\|$
and due to~(\ref{eq.lambdas}) does not change the spectral radius of $X$. 
Then we solve the $m$-dimensional problem 
\begin{equation}\label{eq.v1}
\left\{
\begin{array}{l}
\|X^{(1, 1)} \, - \, A^{(1, 1)}\|\ \to \ \min\\
\rho(X^{(1, 1)})\ \le \ 1 
\end{array}
\right.
\end{equation} 
In other words, we apply the Algorithm to the same problem of a smaller dimension. 
Then we denote its solution by $X^{(1, 1)}$ and set  
\begin{equation}\label{eq.reduc1}
X \quad = \quad \left(
\begin{array}{cc}
X^{(1, 1)} & A^{(1, 2)} \\
\bzero & A^{(2, 2)}  
\end{array}
\right)\ .
\end{equation}
\smallskip 
 
{\tt Case Red-2.} $\ \rho(A^{(2, 2)})\, \ge \,  1$.  In this case we 
apply the  Algorithm independently to the blocks $(1, 1)$ and $(2, 2)$. 
Thus, we solve two independent problems of dimensions $m$ and $d-m$: 
\begin{equation}\label{eq.v2}
\left\{
\begin{array}{l}
\|X^{(1, 1)} \, - \, A^{(1, 1)}\|\ \to \ \min\\
\rho(X^{(1, 1)})\ \le \ 1 
\end{array}
\right.
\ ; \qquad 
\left\{
\begin{array}{l}
\|X^{(2, 2)} \, - \, A^{(2, 2)}\|\ \to \ \min\\
\rho(X^{(2, 2)})\ \le \ 1 
\end{array}
\right.
\end{equation}  
Then we denote the solutions as $X^{(1, 1)}$ and $X^{(2, 2)}$ respectively
and set the matrix~$X$: 
\begin{equation}\label{eq.reduc2}
X  \quad = \quad \left(
\begin{array}{cc}
X^{(1, 1)} & A^{(1, 2)} \\
\bzero & X^{(2, 2)}  
\end{array}
\right)\ .
\end{equation}
This is a stationary point (see Proposition~\ref{p.40} in the next section), and the algorithm terminates. 
\bigskip 

\begin{center}
\textbf{4.4. Illustrative example}
\end{center}
\bigskip

Consider the matrix
\[
A = \left( \begin{array}{ccccc}
    0.7  &  0.2  &  0.1  &  0.5  &  1.0 \\
    0.3  &  0.6  &  0.2  &  0.8  &  0.3 \\
    0.5  &  0.7  &  0.9  &  1.0  &  0.5 \\
    0.1  &  0.1  &  0.3  &  0.8  &  0.3 \\
    0.8  &  0.2  &  0.9  &  0.3  &  0.2
\end{array}
\right), \qquad \mbox{with} \quad \rho(A)=2.4031.
\]

After the first inner optimization step, the following matrix is found:
\[
\widetilde X_1   = \left( \begin{array}{ccccc}
    0.4349  &  0.1406  &  0.0652  &  0.4912  &  0.9345 \\
    0       &  0.3751  &  0.0682  &  0.7668  &  0.0518 \\
    0       &  0.3383  &  0.6881  &  0.9466  &  0.1009 \\
    0       &       0  &       0  &  0.5917  &       0 \\
    0.2989  &  0.0878  &  0.8343  &  0.2834  &  0.0762
\end{array}
\right), 
\]
which is reducible and has distance $1.1894$ from $A$.

A reordering allows to obtain 
\[
X_1 = \left( \begin{array}{ccccc}
    0.0762  &  0.2989  &  0.0878  &  0.8343  &  0.2834 \\
    0.9345  &  0.4349  &  0.1406  &  0.0652  &  0.4912 \\
    0.0518  &       0  &  0.3751  &  0.0682  &  0.7668 \\
    0.1009  &       0  &  0.3383  &  0.6881  &  0.9466 \\
    0       &       0  &       0  &       0  &  0.5917 \\
\end{array}
\right), \quad
A_1 = \left( \begin{array}{ccccc}
    0.2  &  0.8  &  0.2  &  0.9  &  0.3 \\
    1.0  &  0.7  &  0.2  &  0.1  &  0.5 \\
    0.3  &  0.3  &  0.6  &  0.2  &  0.8 \\
    0.5  &  0.5  &  0.7  &  0.9  &  1.0 \\
    0.3  &  0.1  &  0.1  &  0.3  &  0.8 \\
\end{array}
\right).
\]

We are in the first case, {\tt Case Red-1}, that is $\rho(A_1^{(2,2)}) = 0.8 < 1$.
So we continue applying the optimization Algorithm to the matrix $A_1^{(1,1)}$.
This gives the matrix
\[
X_1^{(1,1)} = \left( \begin{array}{cccc}
         0  &  0.4204  &  0.1759  &  0.6770 \\
    0.7343  &  0.3796  &  0.1797  &       0 \\
    0.0274  &       0  &  0.5791  &  0.0069 \\
    0.1334  &  0.0580  &  0.6719  &  0.6403 
		\end{array} \right),		
\]
which is irreducible and primitive. Hence the final matrix is given 
\[
X_2 = \left( \begin{array}{ccccc}
         0  &  0.4204  &  0.1759  &  0.6770 & 0.3 \\
    0.7343  &  0.3796  &  0.1797  &       0 & 0.5 \\
    0.0274  &       0  &  0.5791  &  0.0069 & 0.8 \\
    0.1334  &  0.0580  &  0.6719  &  0.6403 &   1 \\
		     0  &       0  &       0  &       0 & 0.8
\end{array} \right),		
\]
which optimally approximates the matrix $A_1$, that is $A$ expressed in
the new coordinates. The distance is $1.1037$, which is indeed smaller
than $1.1894$. 
In the old coordinates, the computed optimal solution is $X^* = A - \Delta$ is
\[
X^* = \left( \begin{array}{ccccc}
    0.3796  &  0.1797  &       0  &  0.5  &  0.7343 \\
         0  &  0.5791  &  0.0069  &  0.8  &  0.0274 \\
    0.0580  &  0.6719  &  0.6403  &  1.0  &  0.1334 \\
         0  &       0  &       0  &  0.8  &       0 \\
    0.4204  &  0.1759  &  0.6770  &  0.3  &       0
\end{array} \right),
\quad
\Delta = \left( \begin{array}{ccccc}		
    0.3204  &  0.0203  &  0.1000    &     0  &  0.2657 \\
    0.3000  &  0.0209  &  0.1931    &     0  &  0.2726 \\
    0.4420  &  0.0281  &  0.2597    &     0  &  0.3666 \\
    0.1000  &  0.1000  &  0.3000    &     0  &  0.3000 \\
    0.3796  &  0.0241  &  0.2230    &     0  &  0.2000		
\end{array} \right).		
\]


%
%



 %

\bigskip

\begin{center}
\textbf{4.5. Realization and computational costs} 
\end{center}

\bigskip 

An advantage of the algorithm is a relatively low computational cost of each iteration. 
Problem~(\ref{eq.v}) is solved as $d$ separate problems, one in each row. 
In the $i$-th row we find the minimum of $\|\bx_i \, - \, \ba_i\|^2$ under the
constraints $(\bx_i, \bv) \, \le \, v_i, \ \bx_i \ge 0$ ($\ba_i, \bx_i$ are the $i$s rows of the matrices $A$ and $X$ respectively, 
and $\bv = \bv_{k-1}$). 
This is a $d$-dimensional convex quadratic problem and is easily solved by 
quadratic programming. The same problem of minimizing positively definite 
quadratic form on the positive orthant under one linear constrained arises in 
many applications. 
Since the objective function is strictly convex, the solution $\bx_i$ is unique and   
by the Karush-Kuhn-Tukker theorem (see e.g \cite{F}), is characterized by the 
equation: 
\begin{equation}\label{eq.AX}
\bx_i \quad  = \quad  
\left\{
\begin{array}{lcl}
\ba_i \, ,& \mbox{if} & (\ba_i, \bv) \, \le \, v_i\\
\ba_i \, - \, \lambda\, \bv\, + \, \Lambda_i, & \mbox{if} & (\ba_i, \bv) \, > \, v_i\, , 
\end{array}
\right. 
\end{equation}
where $\lambda > 0$ is a multiplier and $\Lambda_i \ge 0$ is a vector orthogonal to $\bx_i$, i.e., 
it has zeros on the positions of positive components of $\bx_i$. If $\bx_i > 0$, then 
$\Lambda_i = 0$ and $\lambda \, = \, \frac{(\ba_i, \bv) \, - \, v_i}{\|\bv\|^2}$. So, in this case 
$\bx_i$ is explicitly computed: 
\begin{equation}\label{eq.AXp}
\bx_i \quad  = \quad  
\left\{
\begin{array}{lcl}
\ba_i & \mbox{if} & (\ba_i, v) \, \le \, v_i\\
\ba_i \ - \  \, \frac{(\ba_i, \bv) \, - \, v_i}{\|\bv\|^2}\ v\, , & \mbox{if} & (\ba_i, v) \, > \, v_i\, . 
\end{array}
\right. 
\end{equation}
In general, if $\bx_i$ has zeros, it is characterized by equation~(\ref{eq.AX}) and is computed numerically.  
For even iterations~$k$, with the problem~(\ref{eq.u}), 
 the formulas are the same, with replacing rows by columns and 
$\bv$ by $\bu$. The most expansive operation is the computing the leading eigenvector of $X_{k-1}$
(left or right one depending on~$k$) in each step. Further conclusions from formulas~(\ref{eq.AX})
are the following. 
\begin{cor}\label{c.20}
For each $k$, there is a vector $\bell_k \in \re^d_+$ and a matrix $\Lambda_{k} \ge 0, \, \Lambda_k \perp 
X_k$ such that 
$X_k \, = \, A\, - \, \bell_k \, \bv^T_{k-1}\, +\, \Lambda_k\, $ if  $k$ is odd, and 
$X_k \, = \, A\, - \, \bu_{k-1}\, \bell_k^T\, \, +\, \Lambda_k\, $ if~$k$ is even. 
\end{cor}
\begin{cor}\label{c.30}
For each $k$, we have $\, {\rm rank}\, (X_k - A - \Lambda_k)\, = \, 1$. 
\end{cor}

If $X_k$ is primitive, then the $k$-th iteration of the algorithm is complete. 
Otherwise $X_k$ has a block cyclic form~(\ref{eq.cyclic}). We 
multiply each block with a positive weight $s_i$ and minimize the distance to $A$ by optimizing those weight.   Thus we get the matrix $X_{k+1}$ and go to the next iteration. 
The problem~(\ref{eq.prim}) of optimizing weights is also easily solvable. We omit the index $k$
and denote by~$X^{(m)}$
the $m$th block of the matrix $X$ and by~$A^{(m)}$ the corresponding 
pattern of the matrix $A$. Problem~(\ref{eq.prim}) becomes 
\begin{equation}\label{eq.prim1}
\left\{
\begin{array}{l}
\sum_{m=1}^r\ \|\, s_m\, X^{(m)} \, - \, A^{(m)}  \, \|^2\quad \to \quad \min\\
s_1\cdots s_r\ \le \ 1 \, .
\end{array}
\right.
\end{equation}
Its solution satisfies the system of Lagrangian  equations 
$$
s_m^2\, \bigl\|X^{(m)}\bigr\|^2 \, - \, s_m\, \bigl(X^{(m)} \, , \, A^{(m)}\, \bigr) \ + \ \lambda \ , 
\quad , \, m = 1, \ldots , r\, .  
$$
This is a union of univariate quadratic equations depending on one parameter $\lambda>0$. 
Each equation has two positive roots. Taking  every time the smallest one as $s_m$, we 
then find the numbers $s_1(\lambda), \ldots , s_r(\lambda)$. Then we find the smallest $\lambda$
for which $s_1\cdots s_r=1$, we find the optimal weights.  

\begin{cor}\label{c.40}
If the matrix $X$ has the optimal weights of the blocks, i.e., for that matrix \linebreak $s_1= \ldots =s_m = 1$, 
then all the scalar products $\bigl(\, X^{(m)}\, , \, X^{(m)} \, - \, A^{(m)}\, \bigr)$ are the same for $m=1, \ldots , r$. 
\end{cor}

\bigskip 

\begin{center}
\textbf{4.6. Optimal stabilization at a stationary point}
\end{center}

\bigskip 

Now we are going to show that if  Algorithm \ref{algo:opt}  stabilizes at some 
matrix $X$, i.e., $X_k = X$ in several subsequent iterations, then $X$ is a  stationary 
point~(Definition~\ref{d.16}). 
In this case the algorithm terminates within finite time. The next step is 
to prove the convergence to a local minimum, this is done in the next subsection.   

Clearly, the distance $\|X_k - A\|$ does not increase in $k$. Moreover,
since each of the problems~(\ref{eq.v}) and~(\ref{eq.u}) possesses a unique solution, it follows that 
  the distance $\|X_k - A\|$ strictly decreases, unless $X_k = X_{k-1}$. If this 
  happens two times in a row, and $X_k$ is primitive, then the algorithm stabilizes at~$X_{k}$.
  If $X_k$ is imprimitive, then one more iteration is needed: finding optimal weights of the blocks, i.e., solving problem~(\ref{eq.prim1}). If it does not change $\|X_k - A\|$, then 
  $X_k$ possesses the optimal weights, that is satisfies conditions of Corollary~\ref{c.40}. 
  Hence, the algorithms stabilizes.  Thus, we have proved the following  
   
\begin{prop}\label{p.20}
The value $\|X_k - A\|$ does not increase in $k$. If it does not change for 
two subsequent iterations at a primitive matrix~$X_k$, or three subsequent 
iterations at an imprimitive one, then the algorithm stabilizes. 
\end{prop}  
 We see that if three consecutive  iterations with the same value of the objective function $\|X-A\|$ mean that the algorithm stabilizes, provided the matrices are irreducible. Now we are going to show that $X$ is a stationary point.  
 
 \begin{theorem}\label{th.10}
If Algorithm \ref{algo:opt} stabilizes at a matrix $X_k$, then~$X_k$ is 
a stationary matrix in the sense of Definition~\ref{d.16}. 
\end{theorem}
To prove the theorem we need some auxiliary results. 
We write $a \sim b$ for two collinear (proportional) vectors.   
\begin{lemma}\label{l.20}
Let $X$ be a  primitive matrix and $a, \tilde a, b, \tilde b$ be non-negative vectors.  If 
the rank-one matrices $a\, b^T$ and $\tilde a\, \tilde b^T$ are equal on ${\rm supp}\, X$, then 
$\tilde a \sim a$ and $\tilde b \sim b$. 
\end{lemma}
The proof is in Appendix. Thus, a rank-one matrix~$C$ has a unique, up to multiplication by a constant, 
presentation $C=ab^T$ on a support of any primitive matrix.

\smallskip 

\begin{prop}\label{p.30}
If in the {\tt Case Red-1} of the algorithm, the matrix $X_{k+1}$
has the same upper triangular form as $X_k$ (with the same sizes and positions of blocks, but with possibly  
new matrices in those blocks), then $X_k$ is stationary.  
\end{prop}
{\tt Proof.} In the $(k+1)$-th iteration of the algorithm we 
compute the left eigenvector $\bu_k$ of the matrix $X_k$ given by formula~(\ref{eq.reduc1})
and solve the problem $\|X - A\|\to \min, \ \bu_k^T X \le \bu_k^T, \ X \ge 0$.
This problem is solved separately in each column:
$\|X^j - A_j\|\to \min \, , \ (\bu_k , \bx^j) \le u_{k, j}$. For $j=d-m+1, \ldots , d$
we already have an optimal solution $\bx^j = \ba^j$, which will not change, because it is unique. 
If the first $m$ columns of $X_{k+1}$ are concentrated in the block $(1, 1)$, then 
$X_{(k+1)}^{(1, 1)}\, = \, X_{k}^{(1, 1)}$, because $X_{k}^{(1, 1)}$ is the solution for this block 
obtained in the previous iteration. Thus, $X_{k} = X_{k+1}$. On the other hand, 
in the next iteration we will have $\bv_{k+1}\, = \, \bv_{k-1}$, and hence 
again $X_{k+2} = X_{k+1}$. Thus, the matrix stays the same for two iterations in a row, hence   
it is stationary.

{\hfill $\Box$}
\smallskip 

\begin{prop}\label{p.40}
If in the case {\tt Case Red-2} of the algorithm, both $X_{k}^{(1, 1)}$
and $X_{k}^{(2, 2)}$ are local minima for their problems, then $X_k$ is a local minimum for the 
original problem~(\ref{eq.prob}). If they both stationary for their problems, 
then $X_k$ is stationary for the 
original problem~(\ref{eq.prob})
\end{prop}
{\tt Proof.} We prove the first part (for the local minima); the proof for stationary matrices  is the same. It suffices to consider the case when the matrices $X_{k}^{(1, 1)}$
and $X_{k}^{(2, 2)}$ are both irreducible. If one of them is reducible, then 
we the same argument to it and the proposition follows by induction in the dimension. 
Adding an arbitrary nonzero matrix $\Delta$ multiplied with small $t > 0$
such that $X_k +t \Delta$ is an admissible matrix. This means 
$X_k +t \Delta \ge 0$ and  $\rho(X_k +t \Delta) \le 1$. 
denote by $\Delta^{(i, j)}$ the corresponding blocks of the matrix $\Delta$. 
Since $\Delta$ is admissible, it follows that $\Delta^{(2, 1)} \ge 0$. 
If $\Delta^{(2, 1)} = 0$, then $\|X_k+ t \Delta - A\| \, \le \, \|X_k - A\|$, 
whenever $t$ is small enough, and hence for  variations $\Delta$ with $\Delta^{(2, 1)} = 0$
the matrix $X_k$ is a local minimum. Indeed, if $\Delta^{(2, 1)} = 0$, then 
the spectral radii of both blocks $(1, 1)$ and $(2, 2)$ of the matrix 
$X_k+ t \Delta$ do not exceed one (because its spectral radius is equal to the maximal 
spectral radius of those two blocks). Hence adding $t \Delta$ with a small $t$ do not reduce 
both $\|X_{k}^{(1, 1)} - A^{(1, 1)} \|$  and $\|X_{k}^{(2, 2)} - A^{(2, 2)}\|$, because 
$X_{k}^{(1, 1)}$
and $X_{k}^{(2, 2)}$ are both  local minima. The value 
$\|X_{k}^{(1, 2)} - A^{(1, 2)} \|$ cannot be reduced either, because this is zero. 
Therefore, it remains to consider the case $\Delta^{(2, 1)} \ne 0$. 
Denote by $\tilde \Delta^{(i, j)}$ the block $(i, j)$
extended by zeros to the whole matrix $\Delta$. 
Since $X^{(1, 1)}$ and $X^{(2, 2)}$ are both irreducible and 
$A^{(1, 2)}\ne 0$ (otherwise $A$ is reducible), we have  
$$
\rho(X_k \, +\, t\, \tilde \Delta^{(2, 1)})\quad  = \quad 1 \ +\ C\, \sqrt{t}\ + \ O(t)\quad  \mbox{as} 
\quad t \, \to \, 0\, ,  
$$
where $C>0$. On the other hand, 
$$
\rho(X_k \, +\, t\,\tilde \Delta^{(i, j)})  \ - \ 1\quad = \quad O\, t \bigl(\|\Delta^{(i, j)}\|\bigr)\, ,   
$$
for every $(i, j)\ne (2, 1)$, because the spectral radius is differentiable 
with respect to $\Delta^{(i, j)}$ for every $(i, j)\ne (2, 1)$. 
Therefore, $\|\Delta^{(i, j)}\| \ge \ C/\sqrt{t}$. 
On the other hand, $A^{(2, 1)} \ne 0$ and hence 
$$
\bigl\|X_k \, +\, t\, \tilde \Delta^{(2, 1)} \, - \, A\, \bigr\|^2 \quad = \quad O(t) \, ,  
$$
while in the other three blocks the square of the distance to $A$ increases at least as 
$C_0\sqrt{t}$, where $C_0$
is a constant. Hence $\bigl\|X_k \, +\, t\, \Delta \, - \, A\bigr\| \, \ge \, 
\bigl\|X_k \, - \, A\bigr\|$.  
 
{\hfill $\Box$}
\smallskip

\smallskip 

{\tt Proof of Theorem~\ref{th.10}}. If $X$ is reducible, the the theorem follows 
by Propositions~\ref{p.30} and~\ref{p.40}. Assume $X$ is irreducible.  Then the eigenvalue~$1$ is simple, 
the leading eigenvectors $u , v$ are well-defined up to multiplication by positive constants, 
and Corollary~\ref{c.20} yields 
\begin{equation}\label{eq.ell}
\begin{array}{lcl}
X & = & A\ - \ \bell_1 \, \bv^T\, +\, \Lambda_1\, \\ 
X & = & A\ - \ \bu\, \bell_{2}^T\, \, +\, \Lambda_{2}
\end{array}
\end{equation}
for some vectors $\bell_i \ge 0$ and matrices $\Lambda_i \ge 0, \, \Lambda_i\, \perp\, X\, \ i = 1,2$. 

Since in the support of the matrix $X$ we have $\Lambda_1 = \Lambda_2 = 0$, 
it follows that $\bell_1 \, \bv^T \, = \,  \bu\, \bell_{2}^T$
on the support of $X$. If $X$ is primitive,  we apply Lemma~\ref{l.20} 
and conclude that $\bell_{1} = r u$ for some $r > 0$. Hence 
$\, X \, = \, A\, - \, r\, \bu\, \bv^T\, +\, \Lambda_1$, which in view of Proposition~\ref{p.10} 
implies that $X \in {\rm locmin}$. This completes the proof 
for primitive $X$.  

If $X$ is imprimitive, then we transfer the matrix $X$ to the cyclic block form~(\ref{eq.cyclic}). 
Respectively, the vector $\bv$ is split into $r$  blocks 
$v = (\bv_1, \ldots , \bv_r)$, where $\bv_i = \bv|_{\Omega_i}$, and the same for 
$\bu = (\bu_1, \ldots , \bu_r)$, 
where $\bu_i = \bu|_{\Omega_i}$. Similarly to the primitive case, we show that 
$A^{(i)}\, - \, X^{(i)}\, = \, \mu_i \, \bu_{i+1} \, \bv_i^T$, on the support of~$X^{(i)}$, where $\mu_i$ are some multipliers, $i = 1, \ldots , r$.
From Corollary~\ref{c.40} it follows that 
$\bigl( X^{(i)}\, , \, \mu_i \bu_{i+1} \, \bv_i^T \bigr)$ is the same
for all~$i$. This scalar product is equal to the trace of the matrix 
$\mu_i\, {X^{(i)}}^T\, , \,  \bu_{i+1} \, \bv_i^T\, = \, \mu_i \bu_{i} \, \bv_i^T \, = \, 
\mu_i\, \bigl(\bu_{i}\, , \, \bv_i \bigr)$. 
We used the fact that $X^{(i)}\bv_i \, = \, \bv_{i+1}$ and $\bu^T_{i+1}X^{(i)}= 
\bu_i^T$.
On the other hand, $\bu^T_{i+1}\bv_{i+1}\, = \, \bu^T_{i+1}X^{(i)}\bv_i\, = \, \bu^T_{i}\bv_{i}$. Thus, all the scalar products 
$(\bu_i, \bv_i), \, i = 1, \ldots , r$ are equal. Therefore, all the numbers $\mu_i$ are equal, hence 
$X \, = \, A \, - \, r \, \bu\, \bv^T \, - \, \Lambda$, and so $\, X\in {\rm locmin}$.

{\hfill $\Box$}
 
\bigskip  
 
 \begin{center}
 \textbf{4.7. Convergence of the algorithm}
 \end{center}
 \bigskip

In the previous subsection we showed that if the Algorithm stabilizes, 
then the point of stabilization is a stationary point. 
Since the value $\|X_k - A\|$ decreases in $k$ and bounded below, 
 it converges as $k \to \infty$. In general, however, it does not imply that the 
 algorithm converges. Theorem~\ref{th.20} below claims that the Algorithm
 indeed converges to a stationary point and, moreover, the rate of convergence is 
at least linear. However, it may not converge to a local minimum. 
For instance, if it starts at a stationary primitive matrix~$X_0$, then 
it stays at $X_0$ forever and stabilizes after the first iteration. 
Say, consider Example~\ref{ex.20}. If $A$ is $2\time 2$ matrix with all entries equal to two and    $X_0 = \frac14 A$. Then $X_0$ is a stationary point and $X_k = X_0$ for all~$k$. Hence, the Algorithm converges to $X_0$, although $X_0$ is not a local minimum. Of course, this situation is not generic and a small variation of $X_0$
may lead to the convergence to a local minumum. In practice, because of roundings, tolerance parameters, etc. such small variations occur in each iterations. Hence, we 
can define the following notion of stable convergence.

 \begin{defi}\label{d.50}
Assume the  Algorithm converges to a matrix~$X$. This convergence is called stable 
(or $X_k$ steadily converges to $X$) if 
there is a number $\varepsilon > 0$ and a number $N\in \n$ such that 
for every $k > N$, the Algorithm starting with a matrix $\tilde X_k$
such that $\|\tilde X_k - X_k\| < \varepsilon$ converges to the same matrix $X$. 
 \end{defi} 
In fact, the stable convergence already implies that the limit point is a local minimum.   
 \medskip 
 
\begin{prop}\label{p.17}
If the convergence is stable, then $X \in {\rm locmin}$. 
\end{prop} 
{\tt Proof}. Assume $X_k$ steadily converges to $X$ and $\|A-X\| = r$, but $X$ is not a local minimum. In this case, we can move $X$ to a distance at most $\varepsilon/2$
so that the distance $\|X-A\|$ decreases by some number $\delta > 0$. This means that 
for all sufficiently big $k$,  we have $\|X_k-A\| < r$. 
Since the convergence is stable, the algorithm starting at $X_k$ has to converge 
to the same limit~$X$. However, this is impossible, because the distance 
to $A$ does not increase each iteration, but finally 
must increase from $\|X_k - A\|$ to $r$.

{\hfill $\Box$}

\medskip

We denote $f(X) = \|X-A\|^2$. The following lemma, whose proof is outlined to Appendix, plays a key role in the proof of convergence. 
\begin{lemma}\label{l.30}
For every $k$, we have $\|X_k-X_{k-1}\|^2 \, \le \, f(X_{k-1}) \, - \, f(X_k)$.  
\end{lemma}

\begin{theorem}\label{th.20}
For arbitrary~$A$, and for an arbitrary choice of the initial matrix~$X_0$, the algorithm converges to a stationary point~$X$ (which may depend on~$X_0$) with the linear rate. This means that there are constants $q\in (0,1)$ and $C>0$
such that $\, \|X_k - X\|\, \le \, C\, q^k, \  k \in \n$.  If the convergence is stable, then $X \in {\rm locmin}$. 
\end{theorem} 
{\tt Proof}. First we show that each limit point of the 
sequence $\{X_k\}_{k\in \z}$ is a point of local minimum. Then 
we prove that this sequence converges to that limit point with a linear rate. 
It suffices to consider the case when the limit point is a primitive matrix, 
the other cases are reduced to this one by the same argument 
as in the proof of Theorem~\ref{th.10}

Applying Lemma~\ref{l.30} and the fact that 
the sequence $f(X_k)$ has a limit  as $k \to \infty$, we see that 
$\|X_k - X_{k-1}\| \to 0$ as $k \to \infty$. 	
By compactness, the sequence $\{X_k\}$ has a limit point $X$.
We assume $X$ is primitive. 
Let $\bv$ and $\bu$ be the right and the left leading eigenvectors of $X$.
They are both strictly positive. For an arbitrary small  $\varepsilon > 0$ and 
for an  arbitrary large $M \in \n$, 
there is a number $m$ such that $\|X_{k-1} - X\| < \varepsilon$
for all $k = m, \ldots , m+2M$. 
Taking $\varepsilon$ small enough, we 
obtain that, on the  support of $X$, the values   
$\|X\, - A\, + \, \bell_{k} \, \bv^T_{k-1}\|\ $ and  $\, 
X \, - \, A\, +\, \bu_{k}\, \bell_{k+1}^T$ are both small and  
$\|\bv_{k-1} - \bv\|$ and $\|\bu_{k-1} -  \bu\|$ are both small 
for all $k = m, \ldots , m+2M$.    
Lemma~\ref{l.20} and the primitivity of $X$ imply that 
$\|\bell_{k-1} - \bu\|$ and $\|\bell_{k+1} - \bv\|$
are both small as well. 
Hence, $\Lambda_k$ has a limit $\Lambda$ as $k \to \infty$, and  $\Lambda \perp X$.  
Thus, $X \, = \, A \ - r\, \bu \, \bv^T$ on the support of $X$. Hence, 
$X \, = \, A \ - \, r\, \bu \, \bv^T \, +\, \Lambda$. Thus, $X$ is a stationary point. 

Now we show that $\|X_k - X\| \, \le \, C\, q^k$ for some $q \in (0,1)$ and $C$. 
 If $\varepsilon$ is smaller than the 
smallest positive entry of $X$, then ${\rm supp}\, X \, \subset \, 
{\rm supp}\, X_k$ for all $k = 1, \ldots , i+N$. 
Denote by  $\bar A$ and $\bar X_k$ the restrictions of those matrices to ${\rm supp}\, X$, i.e., we put all other entries of those matrices 
equal to zero. Similarly, for each $i$, we denote 
$\bar \ba_i$ and $\bar \bv_i$ is the restriction of 
 $\ba_i$ to and $\bv$ to ${\rm supp}\, \bx_i$. The next matrix $X_{k+1}$ is defined from the problem 
 \begin{equation}\label{eq.bar1}
 \left\{
\begin{array}{l}
\|\bar \bx_i \, - \, \bar \ba_i  \|^2 \ \to \ \min\,  \\
(\bar \bx_i\, , \, \bar \bv_i)\, = \,   v_{i} 
\end{array}
\right. 
\end{equation}
 For the solution, we have $\bar \bx_i \, - \, \bar \ba_i \, = \, 
 - \, \bell_i \, \bar \bv_i^T$. The extra Lagrangian  term $\Lambda_i$ 
 vanishes, since $\bar \bx_i$ does not have zeros on the support. 
 Multiplying by $\bar \bv_i$, we get 
 $$
 (\bar \bx_i \, , \, \bar \bv_i ) \ - \ (\bar \ba_i \, , \, 
 \bar \bv_i) \ = \  
 - \, \ell_i \, \|\bar \bv_i\|^2\, ,  
 $$
 where $\, \|\bar \bv_i\|^2\, = \, \sum_{ (\bx_i)_m > 0} v_m^2$. 
Since $(\bar \bx_i \, , \, \bar \bv_i ) \, = \, v_i, \ 
(\bar \ba_i \, , \, \bar \bv_i) \, = \, 
(\bar \ba_i \, , \, \bv)$, we have 
$$
(\bar \ba_i \, , \, \bar \bv ) \ - \ v_i \ = \ 
 - \, \ell_i \, \|\bar \bv_i\|^2\, ,
$$
Therefore, 
 \begin{equation}\label{eq.bar2}
\ell_i \ = \ \frac{(\bar \ba_i \, , \, \bar \bv ) \ - \ v_i}{\|\bar \bv_i\|^2}\ = \ 
\frac{\bigl[(\bar A  \, - \, I)\, \bv\, \bigr]_i }{\|\bar \bv_i\|^2}
\end{equation}
Define the $d\times d$ matrix $B$ as follows: the $i$th row of $B$
is equal to the $i$th row of the matrix $\bar A - I$ divided by 
$\|\bar \bv_i\|^2$. Since, as we have shown above, $\ell = \bu + o(1)$
as $k\to \infty$, equality~(\ref{eq.bar1}) yields 
$\bu_{k+1} \, = \, B \ \bv_{k}\, + \, o(1)$. Similarly, defining 
the matrix $C$: the $j$th row of $C$  
is equal to the $j$th column  of $\bar A - I$ divided by 
$\|\bar \bu_i\|^2$, we obtain $\bv_{k+2} \, = \, C^T \ \bu_{k+1}\, + \, o(1)$.
Iterating we get $\bv_{k+2} \, = \, C^T B\ \bv_{k}\, + \, o(1)$.
Note that $B$ and $C$ are both independent of~$k$. 
Again assuming that $\varepsilon$ is small enough we obtain that the distance between 
 $\bv_{m+2M}$ and $(C^T B)^M\ \bv_{m}$ is small. 
 Taking $m$ and $M$ large enough we see that  $\bv$ is an eigenvector
 of the matrix $C^T B$ corresponding to its eigenvalue $1$ and that 
 all other eigenvalues of this matrix restricted to its eigenspace
 containing all corresponding vectors~$\bv_{k}$ is smaller than one by modulo. 
 If $q$ is the biggest modulus of those eigenvalues, then 
 $q< 1$ and $\|\bv_{k}-\bv\| \, \le \, C\, q^k, \ k \in \n$.
 Arguing similarly for $\bu_k$ and taking into account that the matrix 
 $BC^T$ has the same eigenvalues, we conclude that 
 $\|\bu_{k}-\bu\| \, \le \, C\, q^k, \ k \in \n$. 
 
 Thus, both $\bv_k$ and $\bu_k$ converge to $\bv$ and $\bu$ respectively 
 with the linear rate as $k\to \infty$. Invoking now Corollary~\ref{c.20} and 
 Lemma~\ref{l.20} 
 we see that   $\ell_k \to \bu_k$ with the same rate, and hence 
 $X_k$ converges linearly to $X$.

{\hfill $\Box$}
\medskip

\begin{remark}\label{r.20}
{\em In the proof we see that the rate of linear convergence, i.e., the constant $q$,  is determined by 
the eigenvalues of the matrix $(C^T B)$. If the convergence is stable, then 
the rate is the ratio  between the first and the second largest eigenvalues of
this matrix. When the descent of the function $\|X_k - A\|$ becomes very small, we can compute approximations for $B$ and $C$ and hence, can estimate the~$q$.    
}
\end{remark}
\bigskip 

 \begin{center}
 \textbf{4.8. A favorable case: convergence to a positive matrix}
 \end{center}

\bigskip

Denote $r_k = \|X_k - A\|$. If the $k$th matrix $X_k$  in the Algorithm is strictly positive,
then all formulas are simplified. Assume $k$ is odd (for even $k$
the situation is similar); then    
in  Corollary~\ref{c.20}  we have $\Lambda_k =0$ and  
therefore,    
\begin{equation}\label{eq.xka}
X_k \  = \  A\, - \, \bell_k \, \bv^T_{k-1}
\end{equation}
where $\|\bell_k\| = r_k$. Indeed, 
the eigenvector $\bv_{k-1}$ is normalized to have the unit length, hence
$\|\bu_{k}\, \bell_{k+1}^T \| \, = \, \|\bu_{k}\, \|\, \|\bell_{k+1} \|\, = \, 
\|\bell_{k+1} \|$.  
\begin{prop}\label{p.25} If $X_k > 0$ for some $k$, then 
\begin{equation}\label{eq.xkb}
r_k \ = \ \left\{
\begin{array}{lll}
\|(I - A\ )\bv_{k-1}\| & , & k \ \mbox{\rm is odd}\\
\|(I - A^T)\, \bu_{k}\ \| & , & k \ \mbox{\rm is even}\\
\end{array}
\right.
\end{equation}
\end{prop}
{\tt Proof.} Multiplying both parts of equality~(\ref{eq.xka})
by $\bv_{k-1}$ from the right, we get 
$$
X_k \bv_{k-1}\  = \  A\, \bv_{k-1}\ - \ \bell_k \ (\bv_{k-1}\ , \ 
\bv_{k-1})\ = \ 
A\, \bv_{k-1}\, - \, \bell_k \, .
$$  
On the other hand, $X_k \bv_{k-1} = \bv_{k-1}$. Hence, 
$(I-A)\bv_{k-1} \, = \, - \bell_k$. Consequently,  
$r_k \, = \, \|\bell_k\|\, = \, \|(I-A)\bv_{k-1}\|$. 
For even $k$ the proof is the same. 

{\hfill $\Box$}
\smallskip

\begin{prop}\label{p.50}
If during  $2n+1$ iterations of the Algorithm, all matrices 
$X_0, \ldots , X_{2n}$, are strictly positive, then 
\begin{equation}\label{eq.xk0}
\bv_{2n}  \, = \, a_{2n}\, M^{-n}\bv_{0}\ ; \qquad 
\bu_{2n+1}  \, = \, a_{2n+1}\, N^{-n}\bu_{1}
\end{equation}
where $a_{2n}, a_{2n+1}$ are normalizing constants. 
\end{prop}
{\tt Proof.} Writing~(\ref{eq.xka}) for $k$ and $k+1$ we get the system 
\begin{equation}\label{eq.xk1}
\left\{
\begin{array}{lll}
X_k & = & A\, - \, \bell_k \, \bv^T_{k-1}\,  \\
X_{k+1} & = & A\, - \, \bu_{k}\, \bell_{k+1}^T 
\end{array}
\right.
\end{equation}
Multiplying the first equation by $\bu_{k}^T$ 
from the left and keeping in mind that 
$\bu_k$ is the left leading eigenvector for $X_k$, we 
obtain  
$\, \bu_{k}^T \, = \, \bu_{k}^T A\, \, - \, (\bu_k, \bell_k)\, \bv^T_{k-1}$. 
Therefore 
\begin{equation}\label{eq.xk2}
\mu_k \bv_{k-1} \ = \ (A^T - I)\bu_k\, , 
\end{equation}
where  $\mu_k = (\bu_k, \bell_k)$. Similarly, 
multiplying the second  equation of~(\ref{eq.xk1})
from the right by $\bv_{k+1}$  we 
get  
$\, \bv_{k+1} \, = \,  A\, \bv_{k+1}\, \, - \, \bu_{k}(\bell_{k+1}, \bv_{k+1})$, 
and hence 
\begin{equation}\label{eq.xk3}
\mu_{k+1} \bu_{k} \ = \ (A - I)\bv_{k+1}\, , 
\end{equation}
where  $\mu_{k+1} = (\bell_{k+1}, \bv_{k+1})$. Substituting  $\bu_k$
from~(\ref{eq.xk3}) to~(\ref{eq.xk2}) we obtain 
$\mu_k \mu_{k+1}\, \bv_{k-1}\, = \, (A^T - I)(A - I)\bv_{k+1}$. 
Therefore, $\mu_k \mu_{k+1}\, \bv_{k-1}\, = \, (A^T - I)(A - I)\bv_{k+1}\, = \, 
M\, \bv_{k+1}$. 
Thus, $\bv_{k+1}  \, = \, c_k\, M^{-1}\bv_{k-1}$, where $c_k$ is a constant.
Applying this equality successively for $k=1, 3,  \ldots , 2n-1$ 
we prove the first assertion in~(\ref{eq.xk0}). The second one is 
established in the same way.

{\hfill $\Box$}
\smallskip

\begin{theorem}\label{th.30}
If the Algorithm steadily converges to a strictly positive matrix~$X$, then 
this matrix is a point of global minimum and is explicitly constructed by 
Theorem~\ref{th.7}. Moreover, in this case 
$$
\|X_{n} - X\| \ \le \ C\, \left(\frac{\sigma_1}{\sigma_2}\right)^n\, , 
$$   
where $\sigma_1, \sigma_2$ are the smallest and the second smallest 
singular value respectively of the matrix $I-A$.  
\end{theorem}
{\tt Proof}. If the limit matrix $X$ is strictly positive then 
all the matrices $X_k$ are positive for sufficiently large~$k$. 
Hence, we may assume that the Algorithm starts with a positive matrix~$X_0$ and produces only positive matrices. Proposition~\ref{p.50} implies that 
$\bv_{2n}  \, = \, a_{2n}\, M^{-n}\bv_{0}$. Since the convergence is stable the vector~$\bv_{2n}$ tends to a vector $\bv$, which is an eigenvector of $M^{-1}$
corresponding to its largest eigenvalue, i.e.,  
an eigenvector of $M$ corresponding to its smallest eigenvalue. 
Similarly, $\bu_{2n+1}\to \bu$, where $\bu$ is an eigenvector of $N$ corresponding to its smallest eigenvalue. However, $X_k \to X$, hence $\bu$ and $\bv$ are 
left and right eigenvectors of $X$ respectively. Since each stationary 
point has the form~(\ref{eq.kkt}), the matrix $X$ has the same form with that 
$\bu$ and $\bv$  and with $\Lambda = 0$, because $X>0$.
Thus, $X \, = \, A \, - \, r\, \bu \, \bv^T$.  
We see that all assumptions of Theorem~\ref{th.7} are satisfied, 
hence $X$ is point of global minimum.

It remains to estimate the rate of convergence. 
We have  $\|\bv_{2n} - \bv\| \, \le \, C\, \bigl(\frac{\lambda_2}{\lambda_1}\bigr)^n$, where $\lambda_1, \lambda_2$ are the first and the second largest eigenvalues of $M^{-1}$ respectively. Hence, $\|\bv_{n} - \bv\| \, \le \, C\, \bigl(\frac{\sigma_1}{\sigma_2}\bigr)^n$. The same estimate holds for 
$\bu_n$, and hence for $X_n$ as well.

{\hfill $\Box$}
\smallskip 

\begin{remark}\label{r.30}
{\em The condition of stability of the convergence can not be omitted.
For instance, if $d \ge 3$ and $A = E$ (the matrix of ones), then 
the matrix $X_0 = \frac{1}{d}E$ is a stationary point, although not a 
local minimum  (Example~\ref{ex.20}). 
The Algorithm starting at $X_0$ immediately stabilizes on this matrix, i.e., 
$X_k = X_0$ for all~$k$, hence it converges to~$X_0$.  We see that the Algorithm may 
  converge to a positive matrix which is not a local minimum. Nevertheless, this 
  convergence is unstable, and a small perturbation of the matrix $X_k$ in some   
iteration leads to a different limit. That is why in practice the Algorithm 
converges to a local minimum. This is natural  in view of Proposition~\ref{p.17}.  
}
\end{remark}
\bigskip

\begin{center}
 \textbf{5. How many local minima can occur? }
 \end{center}

\bigskip

Consider an arbitrary ordered partition of the set $\Omega = \{1, \ldots , d\}$
to~$m$ nonempty subsets $\{\Omega_1, \ldots , \Omega_m\}$, where $1\le m\le d$. 
Denote $d_j = |\Omega_j|$.  For an arbitrary $d\times d$   matrix $X$, we denote by 
$X^{(ij)}$ the corresponding $d_i\times d_j$ block in the intersections of rows from $\Omega_i$ and  of columns from $\Omega_j$. 

Let  a  non-negative   matrix $A$ be fixed. For an arbitrary non-negative  matrix $X$, we formulate the following properties: 
\medskip 

\noindent 1) $X$ coincides with $A$
above the diagonal blocks and is zero below them; 
\smallskip 

\noindent 2) for each $j= 1, \ldots , m$, the $j$th diagonal 
block $X^{(j, j)}$  is a stable non-negative  matrix locally closest 
 to $A^{(j, j)}$.  
\smallskip 

\noindent 3) for each $j= 1, \ldots , m$, the $j$th diagonal 
block $X^{(j, j)}$  is the closest 
stable non-negative  matrix to $A^{(j, j)}$. 
\smallskip

In items 2) and 3) closeness is in the set of $d_j \times d_j$ matrices. 
Of course, property 3) is stronger than 2).

Property 1) requires that $X$ is a 
block upper triangular matrix with blocks corresponding to the 
partition of~$\Omega$ and coincides with $A$ above the block diagonal.
So,  the matrix $X$ is uniquely defined out of the diagonal blocks by property 1). 
The diagonal blocks are not defined uniquely even if 3) is satisfied.  
Property 2) implies that  
${\rho(X) = \max\limits_{j=1, \ldots , m}\rho(X^{(j, j)})\, \le \, 1}$.  
\smallskip 

\begin{prop}\label{p.60}
Let a matrix $A$ and an ordered  partition $\{\Omega_i\}_{i=1}^m$ be given. 
If $A$ is strictly positive and $\rho(A^{(j, j)}) > 1$ for all $i = 1, \ldots , m$,  
then every $X$ satisfying 1) and 2) is a locally
closest stable non-negative matrix to~$A$. 
\end{prop}
{\tt Proof}. This proposition  follows from  Proposition~\ref{p.40} directly by applying induction in the number of blocks~$m$.

{\hfill $\Box$}
\smallskip

We call a matrix $A$ {\em lower dominant} if $A_{ij} > A_{ji}$ whenever 
$i > j$. In other words, each component of $A$ below the main diagonal is bigger than 
its reflection above the diagonal. 
\smallskip 
 
\begin{prop}\label{p.70}
Let a matrix $A$ be strictly positive, lower dominant,  and have all its diagonal entries bigger than one.  Then for an arbitrary  ordered  partition 
$\{\Omega_i\}_{i=1}^m$, every matrix $X$ satisfying 1) and 3) is a locally
closest non-negative  stable matrix for~$A$. Moreover, those matrices~$X$ are different for different partitions.  
\end{prop}
{\tt Proof}. Since $a_{ii} > 1$ for all~$i$, it follows that 
the diagonal blocks satisfy~$\rho(A^{(j, j)}) > 1$ for any partition of~$\Omega$.  
Hence, Proposition~\ref{p.60} implies that $X$ is a locally closest stable non-negative matrix to~$A$.  
It remains to show that all those matrices are different for different partitions. 
Assume that the same matrix $X$ corresponds to a different partition
 $\{\Omega'_i\}_{i = 1}^{m'}$. Then either one of the sets $\Omega_j$ is 
 spit by the partition $\{\Omega'_i\}_{i = 1}^{m'}$ into several parts, 
 or one of the sets $\Omega'_j$ is  
 spit by the partition $\{\Omega_i\}_{i = 1}^{m}$. Assume the first case
 (the second one is considered in the same way).  
 In this case, the matrix $X^{(j, j)}$ is block upper triangular, 
 according to the corresponding part of the partition~$\{\Omega'_i\}_{i = 1}^{m'}$
 that splits the set~$\Omega_j$. Thus, the matrix $X^{(j,j)}$
 is zero below the diagonal blocks and coincides with $A^{(j,j)}$
 above them. Denote by $\tilde X^{(j,j)}$ the matrix with the same diagonal blocks 
 as $X^{(j, j)}$ but equal to zero above the diagonal blocks and equal to 
 $A^{(j,j)}$ below them. Clearly, $\rho(\tilde X^{(j, j)}) = \rho(X^{(j, j)}) = 1$. Since $A$ is lower dominant, so is $A^{(j, j)}$, and hence 
 $\|\tilde X^{(j, j)} - A^{(j, j)}\|\, < \, \|X^{(j, j)} - A^{(j, j)}\|$. 
 Hence, $X^{(j, j)}$ is not the closest stable matrix to $A^{(j, j)}$, which contradicts to property 3) in the assumption. 
 
{\hfill $\Box$}
\smallskip

Thus, for every matrix~$A$ satisfying the assumptions of 
Proposition~\ref{p.60}, each ordered partition of the set 
$\{1, \ldots , d\}$ generates its own local minumum of 
the problem~(\ref{eq.prob}) and they are different for different partitions. 
The total number of ordered partitions for a $d$-element set is equal to~$2^d$, 
hence  the problem~(\ref{eq.prob}) has at least $2^d$ different 
points of local minima. Thus, we come to the following conclusion, which justifies
the complexity of the problem:  
 
\begin{cor}\label{c.50}
A strictly positive lower dominant matrix that has all 
diagonal entries bigger than one possesses at least $2^d$ locally 
closest stable non-negative matrices. 
 
\end{cor}
\begin{ex}\label{ex.100}
{\em For a $d\times d$ matrix that has all twos on the diagonal and below it and all ones above the diagonal, problem~(\ref{eq.prob}) has at least $2^d$
local minima. }
\end{ex}
\bigskip  

\begin{center}
\textbf{6. Positive Hurwitz stability}
\end{center}
\bigskip  

All our results can be modified to the Hurwitz stability
of positive systems in a straightforward manner. We will describe the main 
constructions without penetrating the details. 

We recall that a matrix is Hurwitz stable if its  spectral 
abscissa (the maximal real part of eigenvalues) is negative. 
A matrix is called Metzler if all its off-diagonal elements are non-negative. 

Since positive linear systems are defined by Metzler matrices, 
the corresponding problem  are formulated as follows: 
{\em find the closest Hurwitz stable/unstable Metzler matrix to a given matrix~$A$.} 

First of all, the problem of finding the closest stable Metzler 
matrix to a matrix $A$ can be reduced to the case when $A$ 
is Metzler. Otherwise we make the same trick as in the last paragraph 
of the Introduction for non-negative matrices: we define the matrix 
$A_{{\rm Metz}}$ entrywise: on  the diagonal $A_{{\rm Metz}} = A$, and off 
the diagonal $A_{{\rm Metz}} =\max \, \{A, 0\}$. Thus, $A_{{\rm Metz}}$ is a Metzler matrix. Then it is shown easily that the matrices $A$ and $A_{{\rm Metz}}$ 
have the same closest Hurwitz stable Metzler matrix.

The following analogue of the Perron-Frobenius theorem takes place 
for Metzler matrices: the maximal spectral abscissa of a Metzler matrix is 
always attained at a real eigenvalue with a non-negative eigenvector
(leading eigenvector). 
That is why, the Algorithm presented in Section~4 is naturally modified  for computing the closest Hurwitz stable Metzler matrix. In equations 
(\ref{eq.prob-1})  the inequality constraints 
  $X\bv_0 \le \bv_0$ and $X\ge 0$
 are replaced by $X\bv_0 \le 0$  and $x_{i,j} \ge 0, \, i\ne j$, 
 respectively. The same for equations (\ref{eq.prob-2}) and for all iterations
 of the algorithm. 
 All convergence results from Section~4 stay the same and the example 
 from Section~5 is also easily modified for the Hurwitz stability problem.  
 
The closest Hurwitz unstable Metzler matrix is found 
 by the explicit formula $X = A + r\, \bu \bv^T$, where $\bv$ and 
 $\bu$ are the 
 eigenvectors of the matrices $M = A^TA$ and $N = AA^T$ corresponding to their 
 smallest eigenvalues.

\begin{ex}\label{ex.150}
{\em We apply our modified algorithm to the following matrix considered in~\cite{A}:
\[
A = \left( \begin{array}{rrrrr}
    0.6470  &  0.1720  & -0.7490  &  0.7280  &  0.7170 \\
   -0.3540  & -0.0620  & -0.9360  & -0.7730  & -0.7780 \\
    0.0460  &  1.1990  & -1.2690  &  0.8370  &  0.3160 \\
   -0.7930  &  0.8020  &  0.4980  & -1.1280  &  1.4070 \\
   -1.5510  &  1.0530  &  2.7890  & -1.4250  &  0.4010
           \end{array} \right)
\]
The matrix is not Metzler and is unstable since it has $3$ eigenvalue in the right
complex half-plane and spectral abscissa $\alpha(A) \approx 0.5317$. 
Aiming to compute the closest stable Metzler matrix to $A$, 
Anderson found the matrix
\[
X_{\rm Anderson} = \left( \begin{array}{rrrrr}
   -0.0590  &  0.1700  &  0.0030  &  0.6650  &  0.6552 \\
         0  & -0.1730  &  0.0300  &       0  &       0 \\
         0  &  1.1800  & -1.3160  &  0.0080  &       0 \\
         0  &  0.8010  &  0.4950  & -1.1780  &  1.3570 \\
         0  &  1.0400  &  2.7560  &       0  & -0.1830 
\end{array} \right)
\]
whose eigenvalues are all contained in the left complex half-plane and whose spectral abscissa
is $-0.0590$. The distance $\|A - X_{\rm Anderson}\|_F^2 \approx 9.485$.
The algorithm proposed by Anderson makes use of the theory of dissipative Hamiltonian systems, 
which provide a helpful characterization of the feasible set of stable matrices.

Applying our algorithm yields instead the matrix
\[
X^* = \left( \begin{array}{rrrrr}
         0  &  0.1720  &       0  &  0.7280  &  0.7170 \\
         0  & -0.0620  &       0  &       0  &       0 \\
         0  &  1.1990  & -1.3444  &  0.1114  &       0 \\
         0  &  0.8020  &  0.4910  & -1.1956  &  1.2779 \\
         0  &  1.0530  &  2.7535  &       0  & -0.2525
\end{array} \right)
\]
which is quite different from $A_{\rm Anderson}$. 
Its whose eigenvalues are still contained in the left complex half-plane and its  spectral abscissa
is $0$. The distance $\|A - X^*\|_F^2 \approx 9.332$, which slightly improves the 
bound from~\cite{A}.}
\end{ex}

\bigskip 

\begin{ex}\label{ex.160}
{\em 
We next apply our algorithm to the following randomly generated Metzler matrix:
\[
A = \left( \begin{array}{rrrrrr}
    0.5700  &  0.4900  &  0.4700  &  0.7300  &  0.0500  &  0.0200 \\
    0.1400  & -1.1300  &  0.9600  &  0.6700  &  0.3200  &  0.9100 \\
    0.9100  &  0.4500  & -1.7000  &  0.9800  &  0.6000  &  0.1100 \\
    0.8000  &  0.6000  &  0.0400  &       0  &  0.5200  &  0.1400 \\
    0.4800  &  0.5400  &  0.7700  &  0.3600  & -1.0200  &  0.4600 \\
    0.4300  &  0.3300  &  0.9200  &  1.0000  &  0.7600  &  0.0700
           \end{array} \right)
\]
The matrix is very unstable: its 
spectral abscissa is $\alpha(A) \approx 2.1425$.

Applying the first step of our algorithm yields the matrix
\[
X_1 = \left( \begin{array}{rrrrrr}
   -0.4074  &  0       &  0       &  0       &  0       &  0 \\
    0       & -1.2885  &  0.7655  &  0.1879  &  0.0782  &  0.6328 \\
    0.6033  &  0.2551  & -1.9391  &  0.3873  &  0.3028  &  0 \\
    0.4651  &  0.3872  &  0       & -0.6471  &  0.1955  &  0 \\
    0.2254  &  0.3782  &  0.5716  &  0       & -1.2666  &  0.1772 \\
    0       &  0.0376  &  0.5613  &  0.1108  &  0.3141  & -0.4412
\end{array} \right)
\]
whose eigenvalues are all contained in the left complex half-plane and whose spectral abscissa
is $0$. However the matrix is reducible so that we can further optimize it and get
\[
X^* = \left( \begin{array}{rrrrrr}
         0  &       0  &       0  &       0  &       0  &       0 \\
    0.1400  & -1.2885  &  0.7655  &  0.1879  &  0.0782  &  0.6328 \\
    0.9100  &  0.2551  & -1.9391  &  0.3873  &  0.3028  &       0 \\
    0.8000  &  0.3872  &       0  & -0.6471  &  0.1955  &       0 \\
    0.4800  &  0.3782  &  0.5716  &       0  & -1.2666  &  0.1772 \\
    0.4300  &  0.0376  &  0.5613  &  0.1108  &  0.3141  & -0.4412
\end{array} \right)
\]
whose eigenvalues are still contained in the left complex half-plane and whose 
distance from $A$ is improved to $\|A - X^*\|_F^2 \approx 4.690$.}
\end{ex}

\bigskip

\begin{center}
 \textbf{Appendix}
 \end{center}
 
\bigskip

 \noindent {\tt Proof of Lemma~\ref{l.20}.} We call two indices $i, j \in \Omega = \{1, \ldots , d\}$
equivalent if the ratio $a_i : a_j$ is uniquely defined by the equation 
$a\, b^T\, = \, C$ on the set ${\rm supp}\, X$. Thus, the whole set 
$\Omega$ is spit into several equivalence classes $\Omega_1, \ldots , \Omega_r$. 
Denote $X(\Omega_k)\, = \, \{j \ | \ \exists \, i \in \Omega_k \ (i, j)\, \in {\rm supp}\, X \}$. 
Since for each $i \in \Omega_k$, we have $a_jb_i = Y_{ij}, \, (i, j) \in {\rm supp}\, X$, 
it follows that the ratios of all $a_j$  for all $j \in {\rm supp}\, X$ are uniquely 
defined, hence, the ratios of all $a_j \, , \ j \in \Omega_k$ are uniquely defined as well. 
Therefore,  the sets $X(\Omega_k), \, k = 1, \ldots , r$ are just some permutation
of the sets $\Omega_1, \ldots , \Omega_r$. Thus, the matrix $X$ and all its powers 
define permutations of those sets, and hence $X$ cannot be primitive, unless $r=1$. 
Consequently, the ratios of  all entries of the vector $a$ is uniquely defined, 
and hence $\tilde \sim a$. The same with $\tilde b \sim b$.

{\hfill $\Box$}
\smallskip 
 
 {\tt Proof of Lemma~\ref{l.30}.} We  prove this inequality in each row. For every $i = 1, \ldots , d$, we denote by 
$\bx_{i,k-1}$ and $\bx_{i,k}$ the $i$th rows of $X_{k-1}$ and $X_{k}$ respectively. 
We are going to show that   
\begin{equation}\label{obtuse}
\|\bx_{i, k}\, -\, \bx_{i, k-1}\|^2 \ \le \ 
\|\bx_{i, k-1}\, -\, \ba_i\|^2 \ - \ \|\bx_{i, k}\, -\, \ba_i\|^2 
\end{equation}
and then take the sum of those inequalities over $i=1, \ldots , d$. 
Geometrically~(\ref{obtuse}) means that the angle $\angle \, \bx_{i, k} \bx_{i, k-1}\ba_i$ is not acute. Invoking equation~(\ref{eq.AX}) we wee that either 
$\bx_{k, i} \, = \, \ba_{i}$, in which case~(\ref{obtuse}) is obvious, or 
$\bx_{k, i}  \, = \, \ba_i \, - \, \lambda\, \bv_{k-1}\, + \, \Lambda_i$
where $\lambda > 0$ and $\Lambda_i \perp \bx_{k, i}$. 
We have 
$$
\bigl( \ba_i \, - \, \bx_{k, i}\, , \, 
\bx_{k-1, i}\,  -\, \bx_{k, i}\bigr) \ = \ 
\bigl( \, \lambda\, \bv_{k-1}\, - \, \Lambda_i \, , \, 
\bx_{k-1, i}\,  -\, \bx_{k, i}\bigr)\ =\ 
\lambda\, \bigl( \bv_{k-1}\, , \, 
\bx_{k-1, i}\,  -\, \bx_{k, i}\bigr) \ - \  
\bigl( \Lambda_i \, , \, \bx_{k-1, i}\bigr)
$$
Note that $\bigl( \bv_{k-1}\, , \, 
\bx_{k-1, i}\bigr) \, = \, \bigl( \bv_{k-1}\, , \, 
\bx_{k, i}\bigr)\,  = \, (v_{k-1})_i$, therefore 
$$
\bigl( \ba_i \, - \, \bx_{k, i}\, , \, 
\bx_{k-1, i}\,  -\, \bx_{k, i}\bigr) \ = \ 
 \ - \  
\bigl( \Lambda_i \, , \, \bx_{k-1, i}\bigr)\ \le \ 0\, , 
$$
since the vectors $\bx_{k-1, i}$ and $\Lambda_i $ are both non-negative. 
This means that $\angle \, \bx_{i, k} \bx_{i, k-1}\ba_i\, \ge \, 90^o$, which completes the proof. 

{\hfill $\Box$}

\medskip

\section*{Acknowledgments}

Part of this work was developed during some visits to Gran Sasso Science Institute in L'Aquila.
The authors thank the institution for the very kind hospitality.

N. Guglielmi thanks the Italian M.I.U.R. and the INdAM GNCS for financial support and also the Center of Excellence DEWS.


\end{document}